\documentclass[10pt]{article}

\usepackage{a4wide}
\usepackage{amssymb}
\usepackage{amsfonts}
\usepackage{amsmath}
\input xy
\xyoption{arrow} \xyoption{matrix}

\date{}

\newtheorem{proposition}{Proposition}[section]
\newtheorem{theorem}[proposition]{Theorem}
\newtheorem{lemma}[proposition]{Lemma}

\newtheorem{corollary}[proposition]{Corollary}

\def\der{\partial }

\def\nFM0{{\nu }_{F,M_0}}
\def\nFN0{{\nu }_{F,N_0}}
\def\nGN0{{\nu }_{G,N_0}}

\def\N0{ {\bf N}_0 }

\def\t{\otimes}
\def\g{\gamma}
\def\v{\varphi}
\def\ra{\rightarrow}

\def\Xpm{X^{\pm }}

\def\s{\sigma}

\def\l1{{\lambda}_1}

\def\a{\alpha}
\def\a0{ {\alpha }_0}
\def\a1{ {\alpha }_1}

\def\l{\lambda}

\def\nFGM0{{\nu }_{F,G,M_0}}

\def\nFN0{{\nu}_{F,N_0}}

\def\sm{{\sigma}^m}

\def\sm1{{\sigma}^{-1}}

\def\smtp1{{\sigma}^{-t+1}}

\def\S1{S^{-1}}

\def\Xpm1{X^{\pm 1}_1}

\def\sPM1{{\sigma }^{\pm 1}}
\def\sMP1{{\sigma }^{\mp 1 }}

\def\d{\delta}

\def\di{{\rm d.ind}}

\def\L{\Lambda}

\def\G{\Gamma}

\def\Ytm1{Y^{t-1}}
\def\Yim1{Y^{i-1}}

\def\CK{{\cal K}}

\def\CM{{\cal M}}

\def\ass{{\rm ass}}

\def\bQ{\overline{Q}}

\def\ker{ {\rm ker } }

\def\D{ \Delta }

\def\SL2Z{ {\rm SL}_2({\bf Z}) }

\def\CR{ {\cal R}}

\def\Gp1{ G^{1 , 1 } }
\def\P11{ P^{-1 , 1 } }
\def\Pp1{ P^{1 , 1 } }

\def\nCLsr{{}^\nu\kern-2pt {\cal L}^{\sigma , \rho  }}
\def\nP{{}^\nu \kern-2pt P}
\def\nL{{}^\nu\kern-2pt L}
\def\nLL{{}^\nu\kern-2pt \Lambda}
\def\nPsr{{}^\nu\kern-2pt P^{\sigma , \rho  }}
\def\nLsr{{}^\nu\kern-2pt L^{\sigma , \rho  }}
\def\nuCL{{}^\nu\kern-2pt  {\cal L}}
\def\nCLsr{{}^\nu\kern-2pt {\cal L}^{\sigma , \rho  }}
\def\nCL1m{{}^\nu\kern-2pt {\cal L}^{-1 , 1  }}
\def\x1nu{x^\frac{1}{\nu}}
\def\xm1nu{x^{-\frac{1}{\nu}}}

\def\CR{ {\cal R}}

\def\ra{\rightarrow }

\def\CB{{\cal B}}

\def\CC{ {\cal C}}

\def\CP{ {\cal P}}

\def\nAM0{{\nu }_{{\cal A},M_0}}
\def\nAN0{{\nu }_{{\cal A},N_0}}

\def\End{ {\rm End }}

\def\CR{ {\cal R }}
\def\CP{ {\cal P }}

\def\bR{\overline{R}}

\def\bQ{\overline{Q}}

\def\ga{\mathfrak{a}}
\def\gb{\mathfrak{b}}

\def\gn{\mathfrak{n}}
\def\gm{\mathfrak{m}}
\def\gp{\mathfrak{p}}
\def\gq{\mathfrak{q}}

\def\bJ{\overline{J}}

\def\SL{{\rm SL}}

\def\di!{\frac{\der^i}{i!}}
\def\dik!{\frac{\der^k_i}{k!}}

\def\Max{{\rm Max}}

\def\N{\mathbb{N}}

\def\0{\overline{0}}
\def\1{\overline{1}}

\def\Ln1{\L_{n,\overline{1}}}

\def\oa{\overline{a}}

\def\a1{a_{\overline{1}}}

\def\bs{\overline{s}}

\def\S{\Sigma}

\def\vn1{\overrightarrow{n-1}}

\def\im{{\rm im}}

\def\Min{{\rm Min}}

\def\mJ{\mathbb{J}}
\def\mI{\mathbb{I}}

\def\rann{{\rm r.ann}}

\def\bgp{\overline{\gp}}

\def\K1{{\rm K}_1}

\def\hmI1{\widehat{\mI_1}}
\def\tmI1{\widetilde{\mI_1}}
\def\tmJ1{\widetilde{\mJ_1}}
\def\hB1{\widehat{B_1}}
\def\hCB1{\widehat{\CB_1}}

\def\bS{\overline{S}}

\def\Den{{\rm Den}}

\def\Ore{{\rm Ore}}

\def\Den{{\rm Den}}

\def\maxDen{{\rm max.Den}}

\def\br{\overline{r}}
\def\bc{\overline{c}}

\def\bs{\overline{s}}
\def\bt{\overline{t}}
\def\ga{\mathfrak{a}}

\def\udim{{\rm udim}}

\def\gll{\mathfrak{l}}

\def\pCC{{}'\CC}
\def\pCCR{{}'\CC_R}
\def\pCCU{{}'\CC_U}
\def\pCCI{{}'\CC_I}
\def\pCCRJ{{}'\CC (R, J)}
\def\RR{{}_RR}
\def\pQ{{}'Q}
\def\pbQ{{}'\overline{Q}}
\def\pbCC{{}'\overline{\CC}}
\def\pCM{{}'\CM}
\def\bCM{\overline{\CM}}
\def\pSlR{{}'S_l(R)}
\def\pQlR{{}'Q_l(R)}
\def\pS{{}'S}
\def\pga{{}'\ga}
\def\pAsslR{{}'{\rm Ass}_l(R)}
\def\pD{{}'\D}
\def\pG{{}'\G}

\begin{document}

\author{V. V. \  Bavula  
}

\title{The classical left regular left quotient ring of a ring and its semisimplicity criteria}

\maketitle

\begin{abstract}

 Let $R$ be a ring, $\CC_R$ and $\pCCR$ be the set of regular and left regular elements of $R$ ($\CC_R\subseteq \pCCR$). Goldie's Theorem is a semisimplicity  criterion for the classical left quotient ring $Q_{l,cl}(R):=\CC_R^{-1}R$. Semisimplicity  criteria are given for the classical left regular  left quotient ring $'Q_{l,cl}(R):=\pCCR^{-1}R$. As a corollary, two new semisimplicity criteria for $Q_{l,cl}(R)$ are obtained (in the spirit of Goldie).


$${\bf Contents}$$
\begin{enumerate}
\item Introduction.
\item   Preliminaries, proofs of Theorem \ref{28Feb15} and Theorem \ref{B26Feb15}.
 \item Semisimplicity criteria for the ring $\pQ_{l,cl}(R)$.
 \item The left regular left quotient ring of a ring and its semisimplicity criteria.
 \item Properties of $\pSlR$ and $\pQ_{l,cl}(R)$.
 \item The classical left regular left quotient ring of the algebra of polynomial integro-differential operators $\mI_1$.
\end{enumerate}

$\noindent $



\end{abstract}


\section{Introduction}

In this paper, $R$ is a ring with 1,  $R^*$ is its group of units,  module means a left module.

{\bf Semisimplicity criteria for the ring $\pQ_{l, cl} (R)$}.  For each element $r\in R$, let $r\cdot : R\ra R$, $x\mapsto rx$ and $\cdot r : R\ra R$, $x\mapsto xr$. The sets $\pCCR := \{ r\in R\, | \, \ker (\cdot r)=0\}$ and $\CC_R' := \{ r\in R\, | \, \ker (r\cdot )=0\}$ are called the {\em sets of left and right regular elements} of $R$, respectively.  Their intersection $\CC_R=\pCCR \cap \CC_R'$ is the {\em set of regular elements} of $R$. The rings $Q_{l,cl}(R):= \CC_R^{-1}R$ and $Q_{r,cl}(R):= R\CC_R^{-1}$ are called the {\em classical left and right quotient rings} of $R$, respectively. Goldie's Theorem states that the ring $Q_{l, cl}(R)$ is  a semisimple Artinian ring iff the ring $R$ is  semiprime, $\udim (\RR)<\infty$ and the ring $R$ satisfies the a.c.c. on left annihilators ($\udim$ stands for the uniform dimension).

In this paper, we consider/introduce the rings $\pQ_{l, cl} (R) := \pCCR^{-1}R$ (the {\em classical left regular left quotient  ring} of $R$) and  $\pQ_{r, cl} (R) :=R {\CC_R'}^{-1}$ (the {\em classical right regular right quotient  ring} of $R$) and give several semisimplicity criteria for them.  In view of left-right symmetry, it suffices to consider, say `left' case.

A subset $S$ of a ring $R$ is called a {\em multiplicative set} if $1\in S$, $SS\subseteq S$ and $0\not\in S$. Suppose that $S$ and $T$ are multiplicative sets in $R$ such that $S\subseteq T$.   The multiplicative subset $S$ of $T$ is called {\em dense} (or {\em left dense}) in $T$ if for each element $t\in T$ there exists an element $r\in R$ such that $rt\in S$.
Main results of the paper are  semisimplicity criteria for the ring $\pQ_{l, cl}(R)$. For  a left ideal $I$ of $R$, let $\pCC_I:= \{ i\in I\, | \, \cdot i : I\ra I$, $x\mapsto xi$, is an injection$\}$.  For a nonempty subset $S$ of a ring $R$, let $\ass_R(S) := \{ r\in R\, | \, sr=0$  for some $s\in S\}$.

\begin{theorem}\label{28Feb15}
Let $R$ be a ring, $\pCC = \pCC_R$ and $\ga := \ass_R(\pCC )$. The following statements are equivalent.
\begin{enumerate}
\item $\pQ:=\pQ_{l,cl}(R)$ is a semisimple Artinian ring.
\item
\begin{enumerate}
\item $\ga$ is a semiprime ideal of $R$,
\item the set $\pbCC:= \pi (\pCC )$ is a dense subset of $\pCC_{\bR}$
 where $\pi : R\ra \bR := R/ \ga$, $r\mapsto \br := r+\ga$,
\item $\udim ({}_{\bR}\bR ) <\infty$, and
\item $\pCC_V\neq \emptyset$ for all uniform left ideals $V$ of $\bR$.
\end{enumerate}
\item $\ga $ is a semiprime ideal of $R$, $\pbCC$ is a dense subset of $\CC_{\bR}$ and $Q_{l, cl}(\bR )$ is a semisimple Artinian ring.
\end{enumerate}

If one of the equivalent conditions holds then $\pbCC \in \Den_l( \bR , 0)$, $\pbCC$ is a dense subset of $\CC_{\bR }$ and $\pQ \simeq \pbCC^{-1} \bR \simeq Q_{l, cl} (\bR )$. Furthermore, the ring $\pQ$ is a simple ring iff the ideal $\ga$ is a prime ideal.
\end{theorem}

Let $\gn = \gn_R$ be the prime radical of the ring $R$. The following theorem is an instrumental in proving  of several results of the paper including Theorem \ref{28Feb15}. It gives sufficient  conditions for the set $\pCCR$ to be a left denominator set of the ring $R$ such that the ring $\pCCR^{-1}R$ is a semisimple Artinian ring.

\begin{itemize}
\item ({\bf Theorem \ref{A26Feb15}}) {\em Let $R$ be a ring. Suppose that $\udim (\RR ) <\infty$ and $\pCCU\neq \emptyset$ for all uniform left ideals $U$ of $R$. Then
$\pCC_R \in \Den_l(R)$, the ring  $\pQ_{l,cl}(R)$ is a semisimple Artinian ring, $\pQ_{l,cl}(R)\simeq Q_{l, cl} (R/ \ga )$  where $\ga : = \ass_R(\pCC )$ and} $\gn_R\subseteq \ga$.
\end{itemize}

For an arbitrary  ring $R$, the set $\maxDen_l(R)$ of maximal left denominator sets is a non-empty set, \cite{larglquot}. The next semisimplicity criterion for the  ring $\pQ_{l,cl}(R)$  is given via the set $\pCM$ of maximal denominator  sets of $R$ that contain $\pCCR$.

\begin{itemize}
\item ({\bf Theorem \ref{3Mar15}}) {\em Let $R$ be a ring, $\ga = \ass_R(\pCCR )$ and $\pCM :=  \{ S\in \maxDen_l(R)\, | \, \pCCR \subseteq S\}$. The following statements are equivalent.
\begin{enumerate}
\item $\pQ_{l, cl} (R)$  is a semisimple Artinian ring.
\item $\pCM$ is a finite nonempty set, $\bigcap_{S\in \pCM}\ass_R(S)=\ga $, for each $S\in \pCM$, the ring $S^{-1} R$ is a simple Artinian ring and the set $\pbCC :=\{ c+\ga \, | \, c\in \pCCR \}$ is a dense subset of $\CC_{R/ \ga }$ in $R/\ga$.
\end{enumerate}
}
\end{itemize}

Theorem \ref{4Mar15} below is a semisimplicity criterion for the ring $\pQ_{l,cl}(R)$ that is given via the set $\Min_R(\ga )$ of minimal primes of  the ideal $\ga$. Theorem \ref{4Mar15} describes explicitly the set $\pCM$ in Theorem \ref{3Mar15}, see the full version of Theorem \ref{4Mar15} in Section \ref{S3TT3}.

\begin{itemize}
\item ({\bf Theorem \ref{4Mar15}}) {\em Let $R$ be a ring, $\pCC = \pCCR$ and  $\ga = \ass_R(\pCC )$. The following statements are equivalent.

\begin{enumerate}
\item $\pQ_{l, cl} (R)$ is a semisimple Artinian ring.
\item
\begin{enumerate}
\item $\ga$ is semiprime ideal of $R$ and the set $\Min_R (\ga ) $ is a finite set.
\item For each $\gp \in \Min_R(\ga )$, the set $S_\gp := \{ c\in R\, | \, c+\gp \in \CC_{R/\gp}\}$ is a left denominator set of the ring $R$ with $\ass_R(S_\gp ) = \gp$.
\item  For each $\gp \in \Min_R(\ga )$, the ring $S_\gp^{-1}R$ is a simple Artinian ring.
\item The set $\pbCC =\{ c+\ga \, | \, c\in \pCC\}$ is a dense subset of $\CC_{R/\ga }$.
\end{enumerate}
\end{enumerate}
}
\end{itemize}

A ring $R$ is called {\em left Goldie} if it satisfies the a.c.c. on left annihilators and $\udim(\RR )<\infty$.
Theorem \ref{7Mar15} below is a semisimplicity criterion for the ring $\pQ_{l,cl}(R)$  in terms of left Goldie rings.

\begin{itemize}
\item ({\bf Theorem \ref{7Mar15}}) {\em
The following statements are equivalent.
\begin{enumerate}
\item $\pQ_{l,cl} (R)$ is a semisimple Artinian ring.
\item
\begin{enumerate}
\item  $\ga$ is a semiprime ideal of $R$ and the set $\Min_R(\ga )$ is finite.
\item For each $\gp \in \Min_R(\ga )$, the ring $R/ \gp$ is a left Goldie ring.
\item  The set $\pbCC$ is a dense subset of $\CC_{\bR}$.
\end{enumerate}
\end{enumerate}
}
\end{itemize}

Theorem \ref{A7Mar15} is a useful semisimplicity criterion for the ring $\pQ_{l,cl}(R)$ as often we have plenty of simple Artinian localizations of a ring.

\begin{itemize}
\item ({\bf Theorem \ref{A7Mar15}})
{\em  The following statements are equivalent.
\begin{enumerate}
\item $\pQ_{l,cl} (R)$ is a semisimple Artinian ring.
\item There are left denominator sets $S_1, \ldots , S_n$ of the ring $R$ such that
\begin{enumerate}
\item  the rings $S_i^{-1}R$ are simple Artinian rings,
\item  $\ga = \bigcap_{i=1}^n \ass_R(S_i)$, and
\item  $\pbCC$ is a dense subset of $\CC_{\bR}$.
\end{enumerate}
\end{enumerate}
}
\end{itemize}

{\it Remark}. Let $R$ be a ring. If $\pCCR$ is a {\em right} denominator set of the ring $R$  then $\pCCR = \CC_R$ and $Q_{r,cl}'(R)= Q_{r,cl}(R)$ is the classical right quotient ring  of $R$.  Similarly,
if $\CC_R'$ is a {\em left} denominator set of the ring $R$ then $\CC_R' = \CC_R$ and $\pQ_{l,cl}(R)=Q_{l, cl}(R)$ is the classical left quotient ring of $R$.

{\bf Semisimplicity criteria for the  ring $Q_{l,cl}(R)$}.
The next theorem shows that the a.c.c.  condition on {\em left} annihilators in Goldie's Theorem can be replaced by the a.c.c.  condition on {\em right} annihilators (or even by a weaker condition) and adding some extra condition.

\begin{theorem}\label{B26Feb15}
Let $R$ be a ring, $\CC = \CC_R$ and $\pCC = \pCC_R$. The following statements are equivalent.
\begin{enumerate}
\item $Q_{l,cl}(R)$ is a semisimple Artinian ring.
\item $R$ is a semiprime ring, $\udim (\RR )<\infty$, the ring $R$ satisfies the a.c.c. on right annihilators  and $\pCCU\neq \emptyset$ for all uniform left  ideals $U$ of $R$.
\item The ring $R$ is a semiprime ring, $\udim (\RR )<\infty$, the set  $\{ \ker (c_R\cdot ) \, | \, c\in \pCC\}$ satisfies the a.c.c.  and $\pCCU\neq \emptyset$ for all uniform left ideals $U$ of $R$.
\item The ring $R$ is a semiprime ring, $\udim (\RR )<\infty$, the set  $\{ \ker (r_R\cdot ) \, | \, r\in R\}$ satisfies the a.c.c.  and $\pCCU\neq \emptyset$ for all uniform left ideals $U$ of $R$.
    \end{enumerate}
\end{theorem}

Below is another semisimplicity criterion for the ring $Q_{l, cl}(R)$ via $\pCCR$.

 \begin{theorem}\label{A28Feb15}
Let $R$ be a ring. The following statements are equivalent.
\begin{enumerate}
\item $Q_{l, cl} (R)$ is a semisimple Artinian ring.
\item $R$ is a semiprime ring, $\udim (\RR )<\infty$, $\pCCR= \CC_R$ and $\pCCU\neq \emptyset$ for all uniform left ideals $U$ of $R$.
\end{enumerate}
\end{theorem}

Apart from Goldie's Theorem, Theorem \ref{B26Feb15} and Theorem \ref{A28Feb15},  there are several semisimplicity criteria for $Q_{l,cl}(R)$, \cite{Bav-Crit-S-Simp-lQuot}.

{\bf The left regular left quotient ring $\pQ_l(R)$ of a ring $R$ and its semisimplicity criteria}. Let $R$ be a ring. In general, the classical left quotient ring $Q_{l,cl}(R)$ does not exists, i.e. the set of regular elements $\CC_R$ of $R$ is not a left Ore set. The set $\CC_R$ contains the {\em largest left Ore set} denoted by $S_l(R)$ and the ring $Q_l(R) := S_l(R)^{-1}R$ is called the {\em (largest) left quotient ring} of $R$, \cite{larglquot}. Clearly, if $\CC_R$ is a left Ore set then $\CC_R = S_l(R)$ and $Q_{l,cl}(R) = Q_l(R)$. Similarly, the set $\pCCR$ of left regular elements of the ring $R$ is not a left denominator set, in general, and so in this case the classical left regular left quotient ring $\pQ_{l,cl}(R)$ does not exist. The set $\pCCR$ contains the {\em largest} left denominator set $\pSlR$ (Lemma \ref{a8Mar15}.(1)) and the ring $\pQ_l(R):= \pSlR^{-1} R$ is called the {\em left regular left quotient ring} of $R$. If $\pCCR$ is a left denominator set then $\pCCR = \pSlR$ and $\pQ_{l,cl}(R) = \pQ_l(R)$. Theorem \ref{8Mar15} is a semisimplicity criterion for the ring $\pQ_l(R)$.

\begin{itemize}
\item  ({\bf Theorem \ref{8Mar15}})
{\em  Let $R$ be a ring. Then
\begin{enumerate}
\item $\pQlR$ is a left Artinian ring iff $\pQ_{l,cl}(R)$ is a left Artinian ring. If one of the equivalent conditions holds then $\pSlR = \pCCR$ and $\pQlR = \pQ_{l, cl}(R)$.
\item $\pQlR$ is a semisimple  Artinian ring iff $\pQ_{l,cl}(R)$ is a semisimple Artinian ring. If one of the equivalent conditions holds then $\pSlR = \pCCR$ and $\pQlR = \pQ_{l, cl}(R)$.
\end{enumerate}
}
\end{itemize}
So, all the semisimplicity criteria for the ring $\pQ_{l,cl}(R)$ are automatically semisimplicity criteria for the ring $\pQ_l(R)$.

{\bf The rings $\pQ_{l,cl}(\mI_1)$ and $Q_{r,cl}'(\mI_1)$}.
Let $K$ be a field of zero characteristic, $A_n = K\langle x_1, \ldots , x_n$, $\frac{\der}{\der x_1}, \ldots , \frac{\der}{\der x_n}\rangle$ be the {\em Weyl algebra} and $\mI_n = K\langle x_1, \ldots , x_n, \frac{\der}{\der x_1}, \ldots , \frac{\der}{\der x_n}, \int_1,\ldots , \int_n\rangle$ be the {\em algebra of polynomial integro-differential operators}. The ring $Q(A_n):=Q_{l, cl}(A_n)$ is a division ring and $Q_{l, cl}(A_n)=Q_l(A_n)=\pQ_{l, cl}(A_n)=\pQ_l(A_n)$.

\begin{itemize}
\item ({\bf Lemma \ref{b21Mar15}})
{\em 1. For all $K$-algebras $A$ and $n \geq 1$, the rings $Q_{l,cl}(\mI_n\t A)$ and $Q_{r,cl}(\mI_n\t A)$ do not exist.

 2. For all $K$-algebras $A$ and $n \geq 1$, the rings $Q_l(\mI_n\t A)$ are not left Noetherian and  the rings  $Q_r(\mI_n\t A)$  are not right Noetherian.
 }
\end{itemize}

As an application of some of the results of the paper the rings $\pQ_{l,cl}(\mI_1)$ and $Q_{r,cl}'(\mI_1)$ are  found.

\begin{itemize}
\item ({\bf Theorem  \ref{19Mar15}})
{\em  $\pQ_{l,cl}(\mI_1) \simeq Q(A_1)$ and $Q_{r,cl}'(\mI_1) \simeq Q(A_1)$ are division rings.}
\end{itemize}
Explicit descriptions of the sets $\pCC_{\mI_1}$ and $\CC_{\mI_1}'$ are given in Theorem \ref{24Mar15}. This and some other results demonstrate that on many occasions the ring $\pQ_{l,cl}(R)$ has `somewhat better properties' than $Q_{l,cl}(R)$ which for $R=\mI_n$ even does not exist.

{\bf Conjecture}. {\em  $\pQ_{l,cl}(\mI_n) \simeq Q_{l,cl}(A_n)$ is a division ring.}


\section{Preliminaries, proofs of Theorem \ref{28Feb15} and Theorem \ref{B26Feb15}}\label{PRLM}


The following notation is fixed in the paper.

{\bf Notation}:

\begin{itemize}
\item  $R$ is a ring with 1, $\gn =\gn_R$ is its prime radical and $\Min (R)$ is the set of minimal primes of $R$;
\item   $\CC = \CC_R$  is the set of regular elements of the ring $R$ (i.e. $\CC$ is the set of non-zero-divisors of the ring $R$);
\item   $Q_{l,cl}(R):= \CC^{-1}R$ is the {\em classical left quotient ring}  (the {\em classical left ring of fractions}) of the ring $R$ (if it exists) and $Q^*$ is the group of units of $Q$;
  \item $\Ore_l(R):=\{ S\, | \, S$ is a left Ore set in $R\}$ and  $\ass
(S):= \{ r\in R \, | \, sr=0$ for some $s=s(r)\in S\}$;
\item $\Den_l(R):=\{ S\, | \, S$ is a left denominator set in $R\}$;
        \item $\Den_l(R, I )$ is the set of left denominator sets $S$ of $R$ with $\ass (S)=I$ where $I$ is an ideal of $R$;
             \item $\maxDen_l(R)$ is the set of maximal left denominator sets of $R$ (it is always a {\em non-empty} set, \cite{larglquot}).
   \item   $\pCC := \pCCR$  is the set of left regular elements of the ring $R$ and $\ga : = \ass_R(\pCC )$,
\item   $\pQ :=\pQ_{l,cl}(R):= \pCCR^{-1}R$ is the {\em classical left regular left quotient ring} and $\pQ^*$ is the group of units of $Q$;
    \item if $\ga$ is an ideal of $R$ then $\bR := R/\ga$,   $\pi : R\ra \bR$, $r\mapsto \br := r+\ga$, and $\pbCC := \pi (\pCC )$;
        \item $\pSlR$ is the largest left denominator set in $\pCCR$ and $\pQ_l(R):= \pSlR^{-1}R$ is the left regular left quotient ring of $R$;
            \item $\pga := \ass_R(\pSlR )$ and ${}'\pi : R\ra \bR':= R/ \pga$, $r\mapsto \br := r+\pga$.
        \end{itemize}

{\bf Sufficient conditions for semisimplicity of the ring $\pQ_{l,cl}(R)$}.  In this section,   proofs are given of Theorem \ref{28Feb15}, Theorem \ref{B26Feb15} and Theorem \ref{A28Feb15}. Let $I$ be a nonzero left ideal of a ring $R$.  Sufficient conditions are given for a right Noetherian ring to have a semisimple  left quotient ring (Corollary \ref{a26Feb15}). For each ideal $\ga$ of a ring $R$, the left singular ideal $\zeta_l(R, \ga )$ of $R$ over $\ga$ is introduced that, in the case when $\ga =0$, coincides  with the (classical) left singular ideal $\zeta_l(R)$ of $R$. It is proved that $\zeta_l(R, \ga )$ is an ideal of $R$.

 Let $\pCCI := \{ i\in I\, | \, \cdot i : I\ra I$, $x\mapsto xi$ is an injection$\}$.
A nonzero module is called a {\em uniform module} if every two of its nonzero  submodules have nonzero intersection.
\begin{lemma}\label{a27Feb15}
Suppose that $R$ is a ring, $U$ is a left uniform ideal of $R$, $u\in \pCCU$ and $K=\ker (\cdot u_R)$. Then
\begin{enumerate}
\item $U\oplus K$ is an essential left ideal of $R$.
\item If $I$ is a left ideal of $R$ such that $U\subseteq I$ then $U\oplus (K\cap I)$ is an essential left $R$-submodule of $I$.
\end{enumerate}
\end{lemma}

{\it Proof}. 1. Clearly, $U\cap K=0$ (since $\ker(\cdot u_U)=0$ and $Ku=0$). So, $U+K= U\oplus K$. Suppose that the left ideal $J:= U\oplus K$ of $R$ is not essential, we seek a contradiction. Then $J\cap V=0$ for some nonzero left ideal $V$ of $R$. The map $\cdot u_V: V\ra U$ is an injection. So, $Vu\cap Uu\neq 0$, i.e. $vu=u'u$ for some nonzero elements $v\in V$ and $u'\in U$, and so $k:= v-u'\in K$. This means that $0\neq v = u'+k\in V\cap J$, a contradiction.

2. The left ideal $J$ of $R$ is essential (statement 1) and $I\neq 0$. Then the intersection $J\cap I = U\oplus I\cap K$ is an essential left $R$-submodule of $I$.   $\Box $

We say that $\udim (\RR)<\infty$ if there are uniform left ideals $U_1, \ldots , U_n$ of $R$ such that $\oplus_{i=1}^n U_i$ is an essential left ideal of $R$. Then $n=\udim (\RR)$ does not depend on the choice of the uniform left ideals $U_i$ and is called the {\em left uniform dimension} of $R$. Similarly, the {\em right uniform dimension}  $\udim (R_R)$ of $R$ is defined.

Let  $J$ be a nonzero ideal of a ring $R$. Let
 $\pCCRJ := \{ r\in R\, | \, \cdot r : J\ra J$, $x\mapsto xr$ is an injection$\}$. We set $\pCC_0(R, 0):= R$. For $r\in R$, let $\cdot r_R= \cdot r: R\ra R$, $ x\mapsto xr$, and $\cdot r_J: J\ra J$, $y\mapsto yr$.

 The {\em classical left quotient ring} $Q_{l, cl}(R)=\CC_R^{-1}R$ often does not exists, i.e. the set $\CC_R$ is not a left Ore set of $R$.  The set $\CC_R$ contains the {\em largest left Ore set} denoted $S_l(R)$ and the ring $Q_l(R):= S_l(R)^{-1}R$
  is called the {\em (largest) left quotient ring} of $R$, \cite{larglquot}. If $\CC_R\in \Ore_l(R)$ then $\CC_R= S_l(R)$ and $Q_{l,cl}(R) = Q_l(R)$.

  \begin{theorem}\label{5Jul10}
  \cite[Theorem 2.9]{larglquot}
The ring  $Q_l(R)$ is a semisimple ring iff the ring  $Q_{l,cl}(R)$ is a semisimple ring. In this case,  $S_l(R) = \CC_R$ and
$Q_l(R) = Q_{l,cl}(R)$.
\end{theorem}

 The next theorem gives sufficient conditions for the set $\pCCR$  to be a left denominator set of $R$ such that the ring $\pCCR^{-1}R$ is a semisimple Artinian ring.

\begin{theorem}\label{A26Feb15}
Let $R$ be a ring and $\pCC := \pCCR$. Suppose that $\udim (\RR ) <\infty$ and $\pCCU\neq \emptyset$ for all uniform left ideals $U$ of $R$. Then
\begin{enumerate}
\item $\pCC \in \Den_l(R, \ga )$ and $\pQ := \pCC^{-1}R$ is a semisimple Artinian ring (where $\ga : = \ass_R(\pCC )$).
\item $\pQ\simeq Q_{l, cl} (R/ \ga )$, an $R$-isomorphism.
\item Let $\pi =\pi_\ga : R\ra R/ \ga$, $r\mapsto \br = r+\ga$, and $\s : R\ra \pQ$, $r\mapsto \frac{r}{1}$. Then
   \begin{enumerate}
\item $\pCC = \pi^{-1} (\pCC_{R/\ga})\cap \pCC (R, \ga )=\{ c\in R\, | \, \cdot \bc_{R/\ga}$ and $\cdot c_\ga$ are injections$\}$ and $\pCC = \s^{-1} (\pQ^*)\cap \pCC (R, \ga )$.
\item $\pCC_{R/\ga }= \CC_{R/\ga}= R/ \ga \cap \pQ^*$.
\end{enumerate}
\item For all essential left ideals $I$ of $R$, $I\cap \pCC \neq \emptyset$.
\item The prime radical $\gn= \gn_R$ of $R$ is contained in the ideal $\ga$ (In general, $\gn \neq \ga$, eg $R= \mI_1$, $\gn_{\mI_1} =0$ but $\ga = F\neq 0$ is the largest proper ideal of $\mI_1$, Theorem \ref{19Mar15}.(1)).
\end{enumerate}
\end{theorem}

{\it Proof}. 4. We use the following fact repeatedly: {\em Given $R$-modules $K\subseteq L \subseteq M$ such that $K$ is an essential submodule of $L$ and $L$ is an essential submodule of $M$ then $K$ is an essential submodule of $M$}. We also use
 repeatedly Lemma \ref{a27Feb15}.

 Let $I$ be an essential left ideal of $R$. Fix a uniform  left ideal, say $U_1$, of $R$ such that $U_1\subseteq I$ and fix an element $u_1\in \pCC_{U_1}$. Let $K_1:= \ker ( \cdot (u_1)_R)$ and $K_1'= I\cap K_1$. By Lemma \ref{a27Feb15}.(2), $U_1\oplus K_1'$ is an essential left $R$-submodule of $I$, hence $U_1\oplus K_1'$ is an essential left ideal of $R$ such that $K_1'u_1=0$. Repeating the same argument for the left ideal $K_1'$ we will find a uniform $R$-submodule $U_2$ of $K_1'$ and an element $u_2\in \pCC_{U_2}$ such that $U_2\oplus K_2'$ is an essential left $R$-submodule of $K_1'$ where $K_2'=K_1'\cap K_2$ and $K_2:= \ker (\cdot (u_2)_R)$. So, $U_1\oplus U_2\oplus K_2'$ is an essential left ideal of $R$ such that $(U_2\oplus K_2') u_1=K_2'u_2=0$. Repeating the same process several times and using the fact that $n:=\udim (\RR )<\infty$, we will find  uniform submodules $U_1, \ldots , U_n$ of $I$ and elements $u_1\in \pCC_{U_1}, \ldots , u_n\in \pCC_{U_n}$ such that

 (i) $J:=\oplus_{i=1}^n U_i$ is an essential left ideal of $R$, and

 (ii) $J_{i+1}u_i=0$, $i=1, \ldots , n-1$ where $J_s:=\oplus_{t=s}^nU_t$, $s=1, \ldots , n$.

 {\em Claim}: $c=u_1+\cdots +u_n\in \pCC\cap I$.

 Clearly, $c\in I$ since all $u_i\in U_i\subseteq I$. We aim to show that $\ker (\cdot c_R)=0$, i.e. $c\in I\cap \pCC$. Since $J$ is an essential left ideal of $R$ it suffices to show that $\ker (\cdot c_J)=0$ (where $\cdot c_J: J\ra J$, $ x\mapsto xc$). The map $\cdot c_J$ respects the ascending filtration of left ideals
 $$ 0=: J_{n+1} \subset J_n\subset J_{n-1}\subset \cdots \subset J_1=J,$$
 i.e. $J_sc= J_s(u_s+\cdots + u_n) \subseteq J_s$ for all $s$ (by (ii) and since $u_s+\cdots + u_n \in J_s$). Moreover, $J_s/J_{s+1}\simeq U_s$ and the map $\cdot c_{J_s/J_{s+1}}= \cdot (u_s)_{U_s}$ is an injection (since $u_s\in \pCC_{U_s}$). Hence, the map $\cdot c_J$ is an injection. The proof of the Claim and of statement 4 is complete.

1. (i) $\pCC\in \Ore_l(R)$: Given $r\in R$ and $c\in \pCC$, we have to show that $ \pCC r \cap Rc\neq \emptyset$. Since $c\in \pCC$, $\udim (\RR c) = \udim (\RR )$, i.e. $Rc$ is an essential left ideal of $R$. Then the left ideal of $R$, $(Rc:r) := \{ a\in R\, | \, ar\in Rc\}$, is an essential left ideal. By statement 4, we can fix an element $c'\in \pCC \cap (Rc: r)$, and so $c' r= r'c$ for some $r' \in R$.

(ii) $\pCC \in \Den_l(R, \ga )$ {\em where} $\ga := \ass_R(\pCC )$: This follows from the statement (i).

(iii) $\pQ$ {\em is a semisimple Artinian ring}: Since $\udim (\RR ) <\infty$ we can fix a direct sum $J=\oplus_{i=1}^n U_i$ of uniform left ideals of $R$ such that $I:= J\oplus \ga$ is an essential left ideal of $R$. By statement 4, $I\cap \pCC\neq \emptyset$. Hence, $\pQ = \pCC^{-1}I= \oplus_{i=1}^n \pCC^{-1} U_i$ (since $\pCC^{-1} \ga =0$). It suffices to show that each $\pQ$-module $\pCC^{-1} U_i$ is a simple module. Suppose that, say $\pCC^{-1}U_1$, is not a simple $\pQ$-module, we seek a contradiction. Then it contains a proper submodule, say $M$. Since ${}_RU_1$ is essential in ${}_R\pCC^{-1} U_1$, the intersection $U_1'= U_1\cap M$ is a nonzero. The left ideal $I' = U_1'\oplus U_2\oplus \cdots \oplus U_n\oplus \ga$ of $R$ is an essential ideal such that $\pCC^{-1} I'$ is a proper left ideal of $\pQ$ (by the choice of $M$) but, by statement 2, $I'\cap \pCC \neq \emptyset$, and so $\pCC^{-1} I'=\pQ$, a contradiction.

3(a) The second equality in the statement (a) is obvious. Let $\CR$ be the RHS of the second equality in the statement (a). To prove that the first equality holds
 it remains to show that $\pCC = \CR$.

(i) $\pCC \subseteq \CR$: Let $c\in \pCC$. Clearly, $c\in \pCC (R, \ga )$. If $\br \bc =0$ for some elements $\br = r+\ga \in R/ \ga$ then $\br =0$ since $R/ \ga\subseteq \pQ$ and $\bc$ is a unit of the ring $\pQ$. So, $\cdot \bc_{R/ \ga }$ is an injection. Then $c\in \CR$.

(ii) $\pCC \supseteq \CR$: If $c\in \CR$ then $c\in \pCC$ as the map $\cdot c_R$ respects the filtration $0\subseteq \ga \subseteq R$.

The equality $\pCC = \s^{-1} (\pQ^*)\cap \pCC (R, \ga )$ follows from the equality $\pCC = \pi^{-1} (\pCC_{R/ \ga })\cap \pCC (R, \ga )$ and the statement 3(b): $\pCC=\pi^{-1} (\pCC_{R/ \ga })\cap \pCC (R, \ga )=\pi^{-1}(R/ \ga \cap \pQ^*)\cap \pCC (R, \ga )=\s^{-1}(\pQ^*)\cap \pCC (R, \ga )$.

3(b) Since $\pQ$ is a semisimple Artinian ring, $\pCC_{R/ \ga}= \CC_{R/\ga }$ and $\CC_{R/ \ga } = R/ \ga \cap \pQ^*$.

2. By statement 1, $\pi (\pCC ) \in \Den_l(R/ \ga , 0 )$ and $\pi (\pCC )^{-1} (R/ \ga ) \simeq \pQ$ is a semisimple Artinian ring. Hence, $\pi (\pCC ) \subseteq S_l (R)$ and $Q_l(R) \simeq \pQ$ is a semisimple Artinian ring. Then $Q_l(R) \simeq Q_{l, cl} (R)$ is a semisimple Artinian ring (Theorem \ref{5Jul10}).

5. Since $\pQ$ is a left Artinian ring, $\pCC^{-1} I$ is an ideal of the ring $\pQ$ for all ideals $I$ of $R$. Let $\gm$ be a maximal ideal of $\pQ$. Then $\s^{-1} (\gm )$ is a prime ideal of the ring $R$: if $IJ\subseteq \s^{-1}(\gm )$ for some ideals $I$ and $J$ of the ring $R$. We may assume that $\ga \subseteq I$ and $\ga \subseteq J$. Then $\pCC^{-1} I \pCC^{-1}J \subseteq \pCC^{-1}\s^{-1} (\gm ) = \gm$, and so one of the ideals $\pCC^{-1}I$ or $\pCC^{-1}J$ belongs to $\gm$. Say, $\pCC^{-1}I \subseteq \gm$. Then $I \subseteq  \s^{-1}(\pCC^{-1} I)\subseteq \s^{-1} (\gm )$. Now, $\gn \subseteq \cap_{\gm \in \Max (R)}\s^{-1}(\gm )$, and so $\pCC^{-1}\gn  \subseteq \cap_{\gm \in \Max (R)}\pCC^{-1}\s^{-1}(\gm )= \cap_{\gm \in \Max (R)}\gm =0$, i.e. $\gn \subseteq \ga$. $\Box $

Let $X$ be a non-empty subset of a ring $R$. The right ideal $\rann (X):=\{ r\in R\, | \, Xr=0\}$ is called the {\em right annihilator} of $X$. Right ideals of this kind are called {\em right annihilator ideals} of $R$.

$\noindent $

{\bf  Proof of Theorem \ref{B26Feb15}}.  Let $Q= Q_{l, cl}(R)$. Then $R\subseteq Q$.

$(1\Rightarrow 2)$ (i) $\udim (\RR )<\infty$: Since $\RR$ is an essential submodule of ${}_RQ$, $\udim (\RR ) = \udim ({}_RQ)= \udim ({}_QQ)<\infty$.

(ii) {\em The ring $R$ is a semiprime ring}: Let $I$ be an ideal of $R$. Since $Q$ is a left Artinian ring, $\CC^{-1}I$ is an ideal of $Q$. If $I$ is a nilpotent ideal then so is $\CC^{-1}I$ (since $ (\CC^{-1} I)^t = \CC^{-1} I^t$ for all $t\geq 1$). Then $\CC^{-1} I=0$, i.e. $I=0$. Therefore, the ring $R$ is a semiprime ring.

(iii) $\pCCU\neq \emptyset$ {\em for all uniform left ideals $U$ of} $R$: The left ideal $M = QU = \CC^{-1}U $ of $Q$ is a simple left $Q$-submodule of $Q$ and $Q= M\oplus N$ for some left $Q$-submodule $N$ of $Q$. Then $1=e_1+e_2$ where $e_1\in M$ and $e_2\in N$ are idempotents of the ring $Q$. Notice that $M = Qe_1= Me_1$ and $e_1M = e_1Qe_1\simeq \End_Q(M)$ is a division ring. So, for each nonzero element $a\in e_1M$, the nonzero $Q$-module homomorphism $\cdot a: M\ra M$, $x\mapsto xa$, of the simple $Q$-module $M$ is an automorphism. Notice that $ue_1=u$ for all elements $u\in U$. In particular, $Ue_1=U$. The ring $R$ is semiprime and $U\neq 0$. By (ii), $0\neq U^2= (Ue_1)^2= Ue_1Ue_1$, and so $ e_1Ue_1\neq 0$. Fix a nonzero element $a= e_1ue_1= e_1u$ in $e_1Ue_1$ where $u\in U$. The map $\cdot a_M: M\ra M$, $m\mapsto ma$, is a bijection. Hence, the map $\cdot a_U: U\ra U$, $u'\mapsto u'a= u'e_1ue_1= u'u$, is an injection, and so $u\in \pCCU$.


(iv) {\em The ring}  $R$ {\em satisfies the a.c.c. on right annihilators}: Let $X$ be a non-empty subset of $R$. Then $\rann_R(X)= R\cap \rann_Q(X)$. Since $Q$ satisfies the a.c.c. on right annihilators, so does $R$.

$(2\Rightarrow 3)$  This implication is obvious.

$(3\Rightarrow 1)$ In view of Theorem \ref{A26Feb15}.(2), it suffices to show that $\ga =0$. By the assumption,  the set $K:= \{ \ker (c_R\cdot ) \, | \, c\in \pCC \}$ satisfies the a.c.c..  Let $\gb := \ker (c_R\cdot )$ be the maximal element in $K$. Clearly, $\gb \subseteq \ga$. We claim that $\gb = \ga$. Otherwise, fix an element  $a\in \ga \backslash \gb$ and an element $c' \in \pCC$ such that $c' a =0$. Since $\pCC \in \Ore_l(R)$ (Theorem \ref{A26Feb15}.(1)), $c''c'= rc$ for some elements $ c''\in \pCC$ and $r\in R$. Then $c_1:= c''c'\in \pCC$ and $\ker ( c_1\cdot ) \supseteq K + aR\varsupsetneqq K$, a contradiction. Therefore, $\gb = \ga$. The left ideal $\ga c$ is a nilpotent ideal ($\ga c \cdot \ga c=0$). The ring $R$ is a semiprime ring, so $\ga c=0$. Then $\ga =0$ since $c\in \pCC$.

$(2\Rightarrow 4\Rightarrow 3)$ These implications are obvious. $\Box$

The next corollary gives sufficient conditions for a {\em right} Noetherian ring to have a semisimple {\em left} quotient ring.
\begin{corollary}\label{a26Feb15}
Let $R$ be a semiprime, right Noetherian  ring with $\udim (\RR )<\infty$ and $\pCCU\neq \emptyset$ for all uniform left ideals $U$ of $R$. Then $Q_{l , cl} (R)$ is a semisimple Artinian ring.
\end{corollary}

{\it Proof}. The corollary follows from Theorem \ref{B26Feb15}.  $\Box $

\begin{lemma}\label{a7Mar15}
Suppose that $S,T\in \Den_l(R)$ and $S\subseteq T$. Then  the map $\v : S^{-1}R\ra S^{-1}T$, $s^{-1}r\mapsto s^{-1}r$, is a ring homomorphism (where $s\in S$ and $r\in R$).
\begin{enumerate}
\item $\v$  is a monomorphism iff $\ass_R(S)=\ass_R(T)$.
\item $\v$  is a epimorphism iff for each $t\in T$ there exists an element $r\in R$ such that $rt\in S+\ass_R(T)$.
\item $\v$  is a isomorphism iff $\ass_R(S)=\ass_R(T)$ and for each element $t\in T$ there exists an element $r\in R$ such that $rt\in S$.
\item If, in addition, $T\subseteq \pCCR$, then $\v$  is a isomorphism iff $\ass_R(S)=\ass_R(T)$ and for each element $t\in T$ there exists an element $r\in \pCCR$ such that $rt\in S$.
\end{enumerate}
\end{lemma}

{\it Proof}. 1. Statement 1 is obvious.

2. $\v$  is a epimorphism iff for each element $t\in T$, $t^{-1}\in \im (\v )$ iff $t^{-1} =s^{-1}r$ for some elements $s\in S$ and $r\in R$ iff $rt-s\in \ass_R(T)$ iff $rt\in S+\ass_R(T)$.

3. Statement 3 follows from statements 1 and 2 and the fact that an element $a\in R$ belongs to $\ass_R(S)$ iff $sa= 0$ for some element $s\in S$.

4. Statement 4 follows from statement 3 and the inclusions $\ker (\cdot rt)\supseteq \ker (\cdot r)$ and $S\subseteq T\subseteq \pCCR$.  $\Box $

$\noindent $

{\bf Proof of Theorem \ref{28Feb15}}. $(3\Rightarrow 1)$ (i) $\pbCC\in \Den_l(\bR , 0)$: Since $\pbCC\subseteq \CC_{\bR}$, it suffices to show that $\pbCC\in \Ore_l(\bR )$. Given elements $s\in \pbCC$ and $ a\in \bR$. Then $as^{-1} = c^{-1} b\in Q_{l, cl} (\bR )$ for some elements $c\in \CC_{\bR}$ and $b\in \bR$, and so $ca= bs$. The set $\pbCC$ is dense in $\CC_{\bR}$. So, $s_1c\in \pbCC$ for some element $s_1\in \bR$. Now, $s_1c\cdot a = s_1bs$. This means that $ \pbCC\in \Ore_l(\bR )$.

(ii) $\pbCC^{-1} \bR = Q_{l, cl} (\bR )$: The equality follows from Lemma \ref{a7Mar15}.(4)
 in view of (i) and the fact that $\pCC$ is dense in $\CC_{\bR}$.

(iii) $\pCC\in \Ore_l(R)$: Given $s\in \pCC$ and $r\in R$. By (i), $s_1r\equiv r_1s\mod \ga$ for some elements $s_1\in \pCC$ and $r_1\in R$. Since $s_1r- r_1s\in \ga$, we can find an element $s_2\in \pCC$ such that $s_2(s_1r- r_1s)=0$, and so
 $s_2s_1\cdot r=s_2 r_1\cdot s$.

 (iv) $\pCC \in \Den_l(R, \ga )$: By (iii), $\ass_R(\pCC )=\ga$. Since every element of $\pCC$ is left regular, the statement (iv) follows.

 (v) $\pQ\simeq \pbCC^{-1}\bR$ (obvious).

 By (ii) and (v), $\pQ\simeq Q_{l,cl}(\bR)$ is a semisimple Artinian ring.

$(2\Rightarrow 1)$ (i) $\pbCC \cap I\neq \emptyset$ {\em for all essential left ideals $I$ of } $\bR$: By the assumptions (c), (d) and Theorem \ref{A26Feb15}.(4), $\pCC_{\bR}\cap I\neq \emptyset$. Fix an element $c\in \pCC_{\bR}\cap I$. Since $\pbCC$ is dense in $\pCC_{\bR}$, $sc\in \pbCC\cap I$ for some element $s\in \bR$.

(ii) $\pbCC\in \Ore_l(\bR )$: Given elements $a\in \bR$ and $c\in \pbCC$, we have to show that $\pbCC a\cap \bR c\neq \emptyset$. Since $c\in \pbCC$, $\udim ({}_{\bR}\bR c) = \udim ({}_{\bR}\bR )$, the left ideal $\bR c$ is an essential left ideal of $\bR$. Then the left ideal of $\bR$, $(\bR c: a) :=\{ b\in \bR\, | \, ba \in \bR c\}$, is an essential left ideal. By the statement (i),  we can find an element $c' \in \pbCC \cap ( \bR c:a)$, and so $c' a = a' c$ for some element $ a' \in \bR$.

(iii) $\pbCC \in \Den_l(\bR , 0)$: In view of (ii), let us  show that $\ass_{\bR}( \pbCC ) =0$. If $\bc \br =0$ for some elements $ c\in \pCC$ and $r\in R$ (where $\bc = c+\ga \in \pbCC$ and $\br = r+\ga \in \bR$) then $cr \in \ga$. Hence, $c' cr=0$ for some element $c' \in \pCC$, and so $r\in \ga$ and $\br=0$. Therefore, $\ass_{\bR}(\pbCC )=0$. It remains to show that if $\br \bc =0$ for some elements $ c\in \pCC$ and $r\in R$ then $\br =0$, i.e. $r\in \ga$. Clearly, $rc\in \ga$, hence $c_1rc=0$ for some element $c_1\in \pCC$. Since $c\in \pCC$, $c_1r=0$, and so $r\in \ga$, as required.

(iv) $\pCC \in \Ore_l(R)$: Given $s\in \pCC$ and $r\in R$. By (ii), $s_1r\equiv r_1s\mod \ga $ for some elements $ s_1\in \pCC$  and $r_1\in R$. Since $s_1r-r_1s\in \ga$, we can find an element $s_2\in \pCC$ such that $s_2(s_1r-r_1s)=0$, and so $s_2s_1\cdot r=s_2r_1\cdot s$  where $s_2s_1\in \pCC$.

(v) $\pCC \in \Den_l(R, \ga )$: By (iii), $\ass_R(\pCC )= \ga$. Since every element of $\pCC$ is left regular, the statement (v) follows.

(vi) $\pQ\simeq \pbCC^{-1} \bR$ (obvious).

(vii) $\pbQ :=\pbCC^{-1} \bR$ {\em is a semisimple Artinian ring}:  Since $\udim ({}_{\bR}\bR ) <\infty$, we can fix an essential direct sum $I=\oplus_{i=1}^n U_i$ of uniform left ideals $U_i$ of the ring $\bR$. By (i), $I\cap \pbCC\neq \emptyset$. Hence, the $\bR$-module $R/ I$ is $\pbCC$-torsion. Therefore, $\pbQ \simeq \oplus_{i=1}^n \pbCC^{-1} U_i$, an isomorphism of $\pbQ$-modules. It suffices to show that each $\pbQ$-module $V_i = \pbCC^{-1} U_i$ is simple. Suppose that, say $V_1$ is not, we seek a contradiction. Then it contains a proper $\pbQ$-submodule, say $M$. By (iii),  ${}_RU_1$ is an essential submodule of ${}_RV_1$, the intersection $U_1'= U_1\cap M$ is nonzero and $\pbCC^{-1}U_1'=M$. The left ideal $J= U_1'\oplus U_2\oplus \cdots \oplus U_n$ of $R$ is an essential left ideal such that $\pbCC^{-1} J=M\oplus V_2\oplus \cdots \oplus V_n$ is a proper left ideal of the ring $\pbQ$. By (i), $I'\cap \pbCC\neq \emptyset$, and so $\pbCC ^{-1} J= \pbQ$, a contradiction.

$(2\Rightarrow 3)$ We continue the proof of the implication $(2\Rightarrow 1)$.

(viii) $\pbQ= Q_{l, cl}(\bR )$ {\em is a semisimple Artinian ring}: This follows from (iii) and (vii).

(ix) $\pbCC$ {\em is  dense in} $\CC_{\bR}$: In view of the statements (iii) and (viii), this follows from Lemma \ref{a7Mar15}.(3).


$(1\Rightarrow 3)$  $(\alpha )$ $\pbCC \in \Den_l(\bR , 0)$ {\em and} $\pQ\simeq \pbCC^{-1} \bR$ (since $\pCC \in \Den_l(R, \ga )$).

$(\beta )$ $\pQ\simeq Q_{l, cl} (\bR )$ {\em is a semisimple Artinian ring} (by $(\alpha )$ and the simplicity of the ring $\pQ\simeq \pbCC^{-1}\bR$).

$(\g )$ $\ga$ {\em is a semiprime ideal of $R$} (by $(\beta )$).

$(\d )$ $\pbCC$ {\em is dense in } $\CC_{\bR}$: By the statements $(\alpha )$ and $(\beta )$, $\pbCC^{-1} \bR \simeq Q_{l,cl}(\bR)$. Now, the statement $(\d )$ follows from Lemma \ref{a7Mar15}.(3).

$(3\Rightarrow 2)$ Recall that $(1\Leftrightarrow 3)$.

(a) $\ga$ is a semiprime ideal of $R$ (this is given).

(b) $\pbCC$ {\em is dense in} $\pCC_{\bR}$: Repeat the proof of  the above statement $(\d )$ bearing in mind that the statements $(\alpha)$ and $(\beta )$ hold in view of the equivalence  $(1\Leftrightarrow 3)$.


(c) $\udim ({}_{\bR}\bR )<\infty$: This follows from the fact that $Q_{l, cl}(\bR )$ is a semisimple Artinian ring.

(d) $\pCC_V\neq \emptyset$ {\em for all uniform left ideals $V$ of} $\bR$ (by Theorem \ref{B26Feb15}.(2), since $Q_{l, cl}(\bR )$ is a semisimple Artinian ring).  $\Box $

\begin{corollary}\label{a28Feb15}
We keep the notation of Theorem \ref{28Feb15}. If $\pQ$ is a semisimple Artinian ring then $\gn \subseteq \ga$ (where $\gn$ is the prime radical of $R$).
\end{corollary}

{\it Proof}. Repeat the proof of statement 5 of Theorem \ref{A26Feb15}.
 $\Box $

 $\noindent $

{\bf Proof of Theorem \ref{A28Feb15}}. $(1\Rightarrow 2)$ The first three conditions are obvious and the fourth holds by Theorem \ref{B26Feb15}.

$(2\Rightarrow 1)$ This implication follows from Theorem \ref{28Feb15} and we keep its notation. Since $\pCCR=\CC_R$, $\ga =0$ and all the conditions (a) - (d) in Theorem \ref{28Feb15}.(2) hold, the ring  $\pQ =Q_{l, cl} (R)$ is a semisimple Artinian ring, by Theorem \ref{28Feb15}. $\Box $

{\bf The left singular ideal of $R$ over $\ga$}.
For a ring $R$, the set $\zeta_l(R):=\{ r\in R\, | \, Ir=0$ for some essential left ideal $I$ of $R\}$ is called the {\em left singular ideal} of $R$. It is an ideal of $R$.  Let $\ga$ be an ideal of $R$. The set $\zeta_l(R, \ga ):=\{ r\in R\, | \, Ir=0$ for some essential left ideal $I$ of $R$ such that $\ga \subseteq I\}$ is called the {\em left singular ideal of $R$ over} $\ga$. It is an ideal of $R$. Clearly, $\zeta_l(R, 0)= \zeta (R)$ and $\zeta_l(R, \ga )$ is a right ideal of the ring $R$.

\begin{lemma}\label{x2Mar15}
For all ideals $\ga$ of a ring $R$, the right ideal $\zeta_l(R, \ga )$ is an  ideal of the ring $R$.
\end{lemma}

{\it Proof}. Let $r\in \zeta := \zeta_l(R, \ga )$ and $Ir=0$ for some essential left ideal $I$ of $R$ such that $\ga \subseteq I$. Let $r'\in R$. The map $f:=\cdot r' : R\ra R$, $ x\mapsto xr'$, is an $R$-homomorphism. Then $f^{-1} (I):= \{ a\in R\, | \, ar'\in I\}$ is an essential left ideal of $R$ that contains the ideal $\ga$. Moreover, $f^{-1} (I)\cdot r'r\subseteq Ir=0$, and so $r'r\in \zeta$. Therefore, $\zeta$ is an ideal of $R$. $\Box $

\begin{proposition}\label{a2Mar2015}
Let $R$ be a ring such that $\pQ_{l, cl} (R)$ is a semisimple Artinian ring and $\ga := \ass_R(\pCCR )$. Then $\zeta_l(R, \ga ) \subseteq \ga $.
\end{proposition}

{\it Proof}. We keep the notation of Theorem \ref{28Feb15}. Let $r\in \zeta := \zeta_l(R, \ga )$. We have to show that $r\in \ga$. Fix an essential left ideal $I$ of $R$ such that $Ir=0$ and $\ga \subseteq I$. Consider the set $S$ of all left ideals $J$ of $R$ such that $J\subseteq I$ and $\ga \cap J=0$. By Zorn's Lemma, let $J$ be a
 maximal element in $S$. Then the left ideal $\ga + J = \ga \oplus J$ is an essential $R$-submodule of $I$, hence it is also an essential left ideal of $R$ (since $I$ is an essential left ideal of $R$).

 {\em Claim}: $\bJ := \pi (J)$ {\em is an essential left ideal of} $\bR$:

  Suppose that this is not true. Then $\bJ \cap \bR \br'=0$ for some nonzero element $\br'= r'+\ga \in \bR$ where $r'\in R$. The left ideal $\ga \oplus J$ of $R$ is essential. So, $(\ga\oplus J)\cap Rr\neq 0$. Let $r''r'= a+j$ be a nonzero element in the intersection for some elements $r''\in R$, $a\in \ga$ and $j\in J$. Then $ca=0$ for some element $c\in \pCC$, and so $cr''r'=cj\neq 0$ (otherwise, $cj=0$, and so $j\in \ga\cap J=0$, a contradiction). Now, $cj \in J\backslash \{ 0\}$, and so $0\neq \overline{cr''r'}=\overline{cj}\in \bR\br' \cap \bJ=0$, a contradiction. So, $\bJ$ is an essential left ideal of the ring $\bR$. This finishes the proof of the Claim.

  Recall that $Ir=0$ and $J\subseteq I$. In particular, $Jr=0$, and so $\bJ \br =0$ and $\br \in \zeta_l(\bR )$.  By Theorem \ref{28Feb15}, $\pQ :=\pQ_{l, cl}(\bR )$ is a semisimple Artinian ring. Then $\zeta_l(\bR )=0$, by \cite[Theorem 2.3.6]{MR}. Therefore, $\br =0$. This means that $r\in \ga$, as required.  $\Box $


\section{Semisimplicity criteria for the ring $\pQ_{l,cl}(R)$}\label{S3TT3}

In this section, proofs are given of several semisimplicity criteria for the ring $\pQ_{l,cl}(R)$ (Theorem \ref{3Mar15}, Theorem \ref{4Mar15}, Theorem \ref{7Mar15} and Theorem \ref{A7Mar15}). It is shown that the left localization radical $\gll_R$ of a ring $R$ is contained in the ideal $\ga = \ass_R(\pCC_R )$ provided $\pQ_{l,cl}(R)$ is a semisimple Artinian ring (Corollary \ref{a3Mar15}). Theorem \ref{5Mar15} gives sufficient conditions for semisimplicity of $\pQ_{l,cl}(R)$ provided the ring $R/ \gn$  is left Goldie (where $\gn$ is the prime radical of the ring $R$).

 For a ring $R$,  let $\pCC = \pCCR$ and $\pCM :=\maxDen_l(R, \pCC ) := \{ S\in \maxDen_l(R)\, | \, \pCC \subseteq S\}$. The first semisimplicity criterion for $\pQ_{l,cl}(R)$ is given via the set $\pCM$ of maximal left denominator sets of $R$ that contain the set $\pCCR$ of left regular elements of $R$.
\begin{theorem}\label{3Mar15}
Let $R$ be a ring, $\pCC = \pCCR$, $\ga = \ass_R(\pCC )$ and $\pCM =  \{ S\in \maxDen_l(R)\, | \, \pCC \subseteq S\}$. The following statements are equivalent.
\begin{enumerate}
\item $\pQ_{l, cl} (R)$ is a semisimple Artinian ring.
\item $\pCM$ is a finite nonempty set, $\bigcap_{S\in \pCM}\ass_R(S)=\ga $, for each $S\in \pCM$, the ring $S^{-1} R$ is a simple Artinian ring and the set $\pbCC$ is a dense subset of $\CC_{R/ \ga }$ in $R/\ga$.
\end{enumerate}
Let $\bR = R/ \ga$ and $\pi : R\ra \bR$, $r\mapsto \br = r+\ga$. If one of the equivalent conditions holds then

$\noindent $

(a) the map $\pCM \ra \bCM :=\maxDen_l( \bR )$, $ S\mapsto \bS:= \pi (S)$, is a bijection with inverse $T\mapsto T':= \pi^{-1} (T)$.

$\noindent $

(b) For all $S\in \pCM$, $\ga \subseteq \ass_R(S)$ and $\pi (\ass_R(S)) = \ass_{\bR}(\bS)$. For all $T\in \bCM$, $\pi^{-1} (\ass_{\bR}(T))=\ass_R(\pi^{-1}(T))$.
$\noindent $

(c) For all $S\in \pCM$, $S^{-1} R\simeq \bS^{-1} \bR$ is a simple Artinian ring.

$\noindent $

(d) $\pQ_{l,cl}(R) \simeq \prod_{S\in \pCM} S^{-1} R\simeq \prod_{S\in \pCM} \bS^{-1} \bR \simeq \prod_{T\in \bCM} T^{-1} \bR\simeq Q_{l,cl} (\bR )$ (semisimple Artinian rings).

\end{theorem}

{\it Proof}. $(1\Rightarrow (a)-(d))$ To prove the implication, first, we prove statements (i)-(vi) below and from which then  the statements (a)-(d) are deduced.

  By Theorem \ref{28Feb15}, $\ga$ is a semiprime ideal of the ring $R$, $\pbCC\in \Den_l(\bR , 0)$ and $Q_{l, cl} (R) \simeq \pbCC^{-1} \bR \simeq \pQ$ is a semisimple Artinian ring.

(i) $\pbCC \subseteq T$ {\em for all} $T\in \bCM$: Recall that $S_l(\bR )$ is the largest left Ore set of $\bR$ that consists of regular elements of the ring $\bR$. Hence, $$\pbCC \subseteq S_l(\bR )\subseteq T,$$  for all $T\in \bCM$, by  \cite[Proposition 2.10.(1)]{Bav-Crit-S-Simp-lQuot}.

(ii) {\em For all} $T\in \bCM$, $T':=\pi^{-1} (T) \in \Den_l(R)$, $\ga \subseteq \ass_{R}(T')$  {\em and} $T'^{-1} R\simeq T^{-1} \bR$: Since  $\pCC \in \Den_l(R, \ga )$ and $\pCC \subseteq \pi^{-1} (\pbCC ) \subseteq \pi^{-1} (T) = T'$ (by (i)),  the result follows from \cite[Lemma 2.11]{Bav-stronglquot}.

(iii)  {\em Given distinct} $T_1, T_2\in \bCM$,  {\em then} $T_1'\neq T_2'$: Suppose that $T_1'= T_2'$. Then,  $T_1=\pi (T_1') = \pi (T_2') = T_2$, a contradiction.

(iv)  {\em For all} $T\in \bCM$, $T'\in \pCM$: By (ii), $T'\in \Den_l(R)$. Then $T'\subseteq S$ for some $S\in \maxDen_l(R)$. Then $S\in \pCM$, since
$$\pCC \subseteq \pi^{-1} (\pbCC ) \stackrel{{\rm (i)}}{\subseteq}\pi^{-1}(T)= T'\subseteq S.$$ Now, $T=\pi (T') \subseteq \pi (S) = \bS$. Since $S\in \pCM$, we have $\pCC \subseteq S$ and so $\ga \subseteq \ass_R(S)$. Therefore, $\bS\in \Den_l(R, \ass_R(S)/\ga )$, and so $\bS\subseteq T_1$ for some $T_1\in \bCM$. Then $T\subseteq \bS \subseteq T_1$, hence $T= T_1$ (since $T, T_1\in \bCM$) and $T= \bS$. Now, $T' = \pi^{-1} (T) = \pi^{-1} (\bS) \supseteq S \supseteq T'$. Therefore, $T' =S\in \pCM$.

(v)  {\em For all} $S\in \pCM$, $S+\ga = S$  {\em and} $ S = \pi^{-1} (\bS)$: For an arbitrary ring $A$ and its maximal left denominator set $S'$, $S' +\ass_A(S') = S'$, by  \cite[Corollary 2.12]{Bav-stronglquot}. Since $\ga \subseteq \ass_R(S)$, we have $S+\ga = S$, hence $S= \pi^{-1} (\bS )$.

(vi)  {\em For all} $S\in \pCM$, $\bS \in \bCM$: Since $\ga \subseteq \ass_R(S)$, $\bS \in \Den_l( \bR , \ass_R(S)/\ga )$. Therefore, $\bS \subseteq T$ for some $T\in \bCM$. Now, $$S\stackrel{{\rm (v)}}{=}\pi^{-1} (\bS ) \subseteq \pi^{-1} (T) = T' \in \pCM,$$ by (iv). Therefore, $S= T'$, and so $\bS = \overline{T'}=T\in \bCM$.

 Now, we are ready to prove the statements (a) - (d).

 (a) By (iv) and (vi), the maps $\pCM \ra \bCM$, $S\mapsto \bS$, and $\bCM \ra \pCM$, $T\mapsto T' := \pi^{-1} (T)$, are well-defined. They are inverses of one another since $S\ra \bS \ra \pi^{-1} (\bS ) = S$ (by (v)); and $T\ra T' \ra \pi (T' ) = T$ (since $\pi$ is a surjection).

(b) Let $S\in \pCM$. Then $\pCC\subseteq S$, and so $\ga \subseteq \ass_R(S)$. Therefore, $\bS \in \Den_l(\bR , \ass_R(S) / \ga )$, and so $\pi ( \ass_R(S)) = \ass_{\bR} (\bS )$.

If $T\in \bCM$ then $T' \in \pCM$ (by the statement (a)), and so $T= \pi (T' ) \in \Den_l(\bR , \ass_R(T')/\ga )$. It follows that $\pi^{-1} (\ass_{\bR } (T)/\ga )= \ass_R(T')$ (since $\ga \subseteq \ass_R(T'))$.

(c) By the statement (a), for all $S\in \pCM$, $\bS\in \bCM$ and $S^{-1} R\simeq \bS^{-1} \bR$. Since $\bS\in \bCM$, the ring $\bS^{-1} \bR$ is a simple Artinian ring, by  \cite[Theorem 3.1]{Bav-Crit-S-Simp-lQuot} (since $Q_{l,cl}(\bR )\simeq \pQ$ is a semisimple Artinian ring).

(d) \begin{eqnarray*}
 \pQ &\simeq & Q_{l, cl} (\bR ) \;\;\;\;\;\;\;\;\;\;\;\,   {\rm (Theorem\;  \ref{28Feb15})}\\
 &\simeq & \prod_{T\in \bCM} T^{-1} \bR\;\;\; \;\;\;\;  {\rm (\cite[Theorem \; 3.1]{Bav-Crit-S-Simp-lQuot})} \\
 &\simeq & \prod_{S\in \pCM}\bS^{-1} \bR\;\;\;\;\;\;\, {\rm (the \; statement \; (a))} \\
 &\simeq & \prod_{S\in \pCM} S^{-1} R\;\;\;\;\; \;\, {\rm (the \; statement \; (c))}.
\end{eqnarray*}

$(1\Rightarrow 2)$ Recall that  $(1\Rightarrow (a)-(d))$ and statement 2 follows from the statements (a)-(d). In more detail,

(i) $1\leq |\pCM |<\infty$: By the statement (a), $|\pCM |= |\bCM |$ and $\bCM$ is a finite set, by \cite[Theorem 3.1]{Bav-Crit-S-Simp-lQuot}.

(ii) $\bigcap_{S\in \pCM}\ass_R(S)=\ga$: By the statement (b), $\ga \subseteq \ass_R(S)$ for all $S\in \pCM$, and so $\ga \subseteq \ga':= \bigcap_{S\in \pCM}\ass_R(S)$. We have to show that $\ga = \ga'$, that is $\ga'  / \ga =0$. Now,
$$ \ga'  / \ga =\bigcap_{S\in \pCM}\ass_R(S)/\ga =\bigcap_{S\in \pCM}\pi (\ass_R(S))\stackrel{{\rm (b)}}{=}\bigcap_{S\in \pCM}\ass_{\bR} (\bS )\stackrel{{\rm (a)}}{=}\bigcap_{\bS\in \bCM}\ass_{\bR} (\bS )=0, $$
by \cite[Theorem 3.1]{Bav-Crit-S-Simp-lQuot}, since $Q_{l,cl}(\bR )$ is a semisimple Artinian ring (by the statement (d)).

(iii) For each $S\in \pCM$, $S^{-1}R$ is a simple Artinian ring (by the statement (c)).

(iv) The set $\pbCC$ is a dense subset of $\CC_{\bR}$ (by Theorem \ref{28Feb15}).

$(2\Rightarrow 1)$ It suffices to show that the conditions of Theorem \ref{28Feb15}.(3) hold. Since $\ga = \bigcap_{S\in \pCM } \ass_R(S)$, the map
$$ R/\ga \ra \prod_{S\in \pCM}S^{-1}R, \;\; r+\ga \mapsto (\frac{r}{1}, \ldots , \frac{r}{1}), $$ is an injection and the direct product is a semisimple Artinian ring. By \cite[Theorem 6.2]{Bav-Crit-S-Simp-lQuot}, $Q_{l,cl}(R/\ga )$ is a semisimple Artinian ring. In particular, the ideal $\ga $ is a semiprime ideal of $R$. Now, by Theorem \ref{28Feb15}.(3), $\pQ$ is a semisimple Artinian ring.  $\Box $

For a ring $R$, the ideal $\gll_R := \bigcap_{S\in\maxDen_l(R)}\ass_R(S)$ is called the {\em left localization radical} of $R$, \cite{larglquot}.
\begin{corollary}\label{a3Mar15}
Let $R$ be a ring such that $\pQ_{l,cl}(R)$ is a semisimple Artinian ring. Then $\gll_R\subseteq \ga$ (where $\gll_R$ is the left localization radical of $R$ and $\ga = \ass_R(\pCCR)$).
\end{corollary}

{\it Proof}. $\gll_R = \bigcap_{S\in\maxDen_l(R)}\ass_R(S) \subseteq \bigcap_{S\in\pCM}\ass_R(S)=\ga$, by Theorem \ref{3Mar15}. $\Box $

Let $R$ be a ring and $I$ be an ideal of $R$. We denote by  $\Min_R(I)$
 the set of minimal prime ideals over $I$. The map $\Min_R(I)\ra \Min (R/I)$, $\gp\mapsto \gp / I$, is a bijection with the inverse $\gq \mapsto \pi_I^{-1}(\gq )$ where $ \pi_I: R\ra R/I$, $r\mapsto r+I$.

The second semisimplicity criterion for $\pQ_{l,cl}(R)$ is given via the minimal primes over $\ga = \ass_R(\pCCR )$. It also gives an explicit description of the elements of the set $\pCM$ (see Theorem \ref{3Mar15} for a definition of $\pCM$).
\begin{theorem}\label{4Mar15}
Let $R$ be a ring, $\pCC = \pCCR$ and  $\ga = \ass_R(\pCC )$. We keep the notation of Theorem \ref{3Mar15}.

\begin{enumerate}
\item $\pQ_{l, cl} (R)$ is a semisimple Artinian ring.
\item
\begin{enumerate}
\item $\ga$ is semiprime ideal of $R$ and the set $\Min_R (\ga ) $ is a finite set.
\item For each $\gp \in \Min_R(\ga )$, the set $S_\gp := \{ c\in R\, | \, c+\gp \in \CC_{R/\gp}\}$ is a left denominator set of the ring $R$ with $\ass_R(S_\gp ) = \gp$.
\item  For each $\gp \in \Min_R(\ga )$, the ring $S_\gp^{-1}R$ is a simple Artinian ring.
\item The set $\pbCC :=\{ c+\ga \, | \, c\in \pCC\}$ is a dense subset of $\CC_{R/\ga }$.
\end{enumerate}

\end{enumerate}
Let $\bR = R/ \ga$ and $\pi : R\ra \bR$, $r\mapsto \br = r+\ga$. If one of the equivalent conditions holds then

$\noindent $

(i) $\pCM =\{ S_\gp\, | \, \gp\in \Min_R(\ga )\} $ and $\ass_R(S_\gp ) = \gp$.

$\noindent $

(ii) $\bCM =\{ S_{\bgp} \; | \, \bgp \in \Min ( \bR )\}$  where   $S_{\bgp}:=\{ \br \in \bR \, | \, \br + \bgp \in \CC_{\bR / \bgp}\}$ and  $\ass_{\bR} (S_{\bgp})=\bgp $.
$\noindent $

(iii) For all $\gp \in \Min_R(\ga )$, $S_\gp^{-1} R\simeq S_{\bgp}^{-1} \bR$ is a simple Artinian ring.

$\noindent $

(iv) $\pQ_{l,cl}(R) \simeq \prod_{\gp\in \Min_R(\ga )} S_\gp^{-1} R\simeq \prod_{\gp\in \Min (\bR )} S_{\bgp}^{-1} \bR \simeq Q_{l,cl} (\bR )$ (semisimple Artinian rings).

\end{theorem}

{\it Proof}. $(1\Rightarrow 2)$. By the assumption, $\pQ$ is a semisimple Artinian ring. By Theorem \ref{28Feb15}, $\ga$ is semiprime ideal of $R$, the set $\pbCC$ is dense in $\CC_{\bR}$ (this is the condition (d)) and the ring $Q_{l, cl}(\bR )$ is a semisimple Artinian ring.

By \cite[Theorem 4.1]{Bav-Crit-S-Simp-lQuot}, $\Min (\bR )$ is a finite set, $\bCM = \{ S_{\bgp} \, | \, \bgp\in \Min (\bR )\}$, $\ass_{\bR}(S_{\bgp})= \bgp$ and $S_{\bgp}^{-1} \bR $ is a simple Artinian ring  for all $\bgp\in \Min (\bR )$.  By Theorem \ref{3Mar15}, $|\Min (\ga ) |= |\Min (\bR )|<\infty$, $\pCM =\{ S_\gp\, | \, \gp\in \Min_R(\ga )\} $,  $\ass_R(S_\gp ) = \gp$ and $S_\gp^{-1}R\simeq S_{\bgp}^{-1}\bR$ is a simple Artinian ring for all $\gp \in \Min_R(\ga )$. Therefore, the properties (a)-(d) hold.

$(2\Rightarrow 1)$  It suffices to show that the conditions of statement 3 of Theorem \ref{28Feb15} hold.

 By the statement (a), the  ideal $\ga$ is a semiprime ideal of $R$. By the statements  (a) and (c), the map
 $$R/\ga \ra \prod_{\gp \in \Min (R, \ga )}S_\gp^{-1}R, \;\; r+\ga \mapsto (\frac{r}{1},\ldots , \frac{r}{1}), $$
is a ring monomorphism. The direct product above  is a semisimple Artinian ring, by the statement (c). By  \cite[Theorem 6.2]{Bav-Crit-S-Simp-lQuot}, the ring $Q_{l,cl}(\bR )$ is a semisimple Artinian ring. By Theorem \ref{28Feb15}, $\pQ$ is  a semisimple Artinian ring. So, the implication $(2\Rightarrow 1)$ holds.

Now, the statements (i)-(iv) follow from Theorem \ref{3Mar15}. $\Box$

A ring $R$ is called {\em left Goldie} if it satisfies the a.c.c. on left annihilators and $\udim(\RR )<\infty$. The third semisimplicity criterion for $\pQ_{l,cl}(R)$ reveals its `local nature' and is given via the rings $R/\gp$ where $\gp\in \Min (\ga )$.

\begin{theorem}\label{7Mar15}
We keep the notation of Theorem \ref{4Mar15}. The following statements are equivalent.
\begin{enumerate}
\item $\pQ_{l,cl} (R)$ is a semisimple Artinian ring.
\item
\begin{enumerate}
\item  $\ga$ is a semiprime ideal of $R$ and the set $\Min_R(\ga )$ is finite.
\item For each $\gp \in \Min_R(\ga )$, the ring $R/ \gp$ is a left Goldie ring.
\item  The set $\pbCC$ is a dense subset of $\CC_{\bR}$.
\end{enumerate}
\end{enumerate}
\end{theorem}

{\it Proof}. $(1\Rightarrow 2)$ Suppose that the ring $\pQ = \pQ_{l, cl}(R)$ is a semisimple Artinian ring. By Theorem \ref{4Mar15}, the conditions (a) and (c) hold, and for each $\gp \in \Min_R(\ga )$, the rings $S_\gp^{-1}R$ are simple Artinian rings (the statement (iii) of Theorem \ref{4Mar15}). Let $\pi_\gp : R\ra R/ \gp$, $r\mapsto r+\gp$. Then $\pi_\gp (S_\gp)\in \Den_l(R/ \gp , 0)$ (since $\ass_R(S_\gp )=\gp$, the statement (i) of Theorem \ref{4Mar15}) and $\pi_\gp (S_\gp)^{-1} (R/ \gp )\simeq S_\gp^{-1} R$ is a simple Artinian ring. Then, $Q_{l,cl}(R/ \gp )\simeq \pi_\gp (S_\gp )^{-1} (R/ \gp )$ is a simple Artinian ring. So, the statement (b) holds.

 $(2\Rightarrow 1)$ Suppose that the conditions (a)-(c) hold. The conditions (a) and (b) means that the ring $\bR = R/ \ga$ is a semiprime ring with $|\Min (\bR )|= |\Min_R(\ga )|<\infty$. By \cite[Theorem 5.1]{Bav-Crit-S-Simp-lQuot}, $Q_{l,cl}(R/ \ga )$ is a semisimple Artinian ring. By Theorem \ref{28Feb15}, $\pQ$ is a semisimple Artinian ring.  $\Box $

The fourth semisimplicity criterion for $\pQ_{l,cl}(R)$ is useful in applications as usually there are plenty of `nice' left denominator sets.

\begin{theorem}\label{A7Mar15}
We keep the notation of Theorem \ref{7Mar15}. The following statements are equivalent.
\begin{enumerate}
\item $\pQ_{l,cl} (R)$ is a semisimple Artinian ring.
\item There are left denominator sets $S_1, \ldots , S_n$ of the ring $R$ such that
\begin{enumerate}
\item  the rings $S_i^{-1}R$ are simple Artinian rings,
\item  $\ga = \bigcap_{i=1}^n \ass_R(S_i)$, and
\item  $\pbCC$ is a dense subset of $\CC_{\bR}$.
\end{enumerate}
\end{enumerate}
\end{theorem}

{\it Proof}. $(1\Rightarrow 2)$ By Theorem \ref{3Mar15}, it suffices to take $\pCM = \{ S_1, \ldots , S_n\}$.

 $(2\Rightarrow 1)$ Suppose that the conditions (a)-(c) hold. By the conditions (a) and (b),  the map
 $$R/\ga \ra \prod_{i=1}^nS_i^{-1}R, \;\; r+\ga \mapsto (\frac{r}{1},\ldots , \frac{r}{1}), $$
is a ring monomorphism. The direct product above  is a semisimple Artinian ring, by the statement (c). By  \cite[Theorem 6.2]{Bav-Crit-S-Simp-lQuot}, the ring $Q_{l,cl}(\bR )$ is a semisimple Artinian ring. So, the conditions of Theorem \ref{28Feb15}.(3) hold, and so  $\pQ_{l, cl}(R)$ is  a semisimple Artinian ring, by Theorem \ref{28Feb15}.  $\Box $

{\bf Sufficient conditions for semisimplicity of $\pQ_{l,cl}(R)$ when $R/ \gn$ is a left Goldie ring}.
Let $R$ be a ring and $I$ be its ideal. Let $\Min (R, I):=\{ \gp \in \Min (R)\, | \, \gp\supseteq I\}$. An important case  for applications is the one when the ring $R/\gn$ is a left Goldie ring, and therefore $Q_{l,cl} (R/ \gn )$ is a semisimple Artinian ring. In this case, the next theorem gives sufficient conditions for semisimplicity of the ring $\pQ_{l,cl}(R)$.

 \begin{theorem}\label{5Mar15}
Let $R$ be a ring, $\pCC = \pCCR$ and  $\ga = \ass_R(\pCC )$. Suppose that the ring $Q_{l, cl}(R/ \gn )$ is a semisimple Artinian ring such that $\ga =\bigcap_{i=1}^n\gp_i$ for some $\gp_1, \ldots , \gp_n\in \Min (R)$.
 If the set $\pbCC :=\{ c+\ga \, | \, c\in \pCCR \}$ is dense in $\CC_{R/ \ga }$ then $\pQ_{l,cl}(R)$ is a semisimple Artinian ring, $\Min_R(\ga ) = \{ \gp_1, \ldots ,\gp_n\}$ and $Q_{l, cl}(R/ \gn )\simeq \prod_{i=1}^n S_{\gp_i}^{-1}R$ where $S_{\gp_i}^{-1}R$ are simple Artinian rings and $S_{\gp_i}:=\{ c\in R\, | \, c+\gp_i\in \CC_{R/ \gp_i}\}\in \Den_l(R,\gp_i)$.
\end{theorem}

{\it Proof}. It is obvious that $\Min_R(\ga )=\{ \gp_1, \ldots , \gp_n\}$. The ring $Q_{l, cl}(R/ \gn )$ is a semisimple Artinian ring. By \cite[Theorem 4.1]{Bav-Crit-S-Simp-lQuot}, $\Min(R/\gn )=\{ \bgp_1,\ldots , \bgp_m \}$ for some $m\geq n$ where $\bgp_i = \gp_i / \gn$ and $\Min (R) = \{ \gp_1, \ldots , \gp_m\}$, $Q_{l,cl}(R/ \gn )\simeq \prod_{i=1}^m S_{\bgp_i}^{-1} (R/ \gn )$ where $S_{\bgp_i}^{-1} (R/ \gn )$ are simple Artinian rings and $S_{\bgp_i}:= \{ c\in R/ \gn \, | \, c+\gp_i/\gn\in \CC_{R/ \gp_i}\}$. The map
$$ R/ \ga \simeq (R/ \gn)/\bigcap_{i=1}^n\bgp_i\ra \prod_{i=1}^n S_{\bgp_i}^{-1} (R/ \gn ), \;\; r+\ga \mapsto (\frac{r}{1}, \ldots , \frac{r}{1}),$$
is a monomorphism and the direct product is a semisimple Artinian ring. Since $\ga / \gn = \bigcap_{i=1}^n \bgp_i$, the conditions (a)-(c) of Theorem \ref{A7Mar15} hold (where $S_i=S_{\bgp_i}$), and so $\pQ_{l,cl}(R)$ is \ semisimple Artinian ring (by Theorem \ref{A7Mar15}). The rest follows from Theorem \ref{4Mar15}. $\Box $


\section{The left regular left quotient ring of a ring and its semisimplicity criteria}\label{LRLQRSSC}

The aim of this section is to prove Theorem \ref{8Mar15} and to establish some relations between the rings $Q_l(R)$ and $\pQ_l(R)$ (Lemma \ref{b8Mar15}, Proposition \ref{1Mar15} and Corollary \ref{a9Mar15}). In particular, to show that  the rings $Q_l(R)$ and $\pQ_l(R)$ are $R$-isomorphic iff $S_l(R) = \pSlR$ (Proposition \ref{1Mar15}.(4)). At the end of the section,  some applications are given for the algebras of polynomial integro-differential operators.

{\bf  The left regular left quotient ring $\pQ_l(R)$ of a ring $R$}. Let $R$ be a ring. Its {\em opposite ring} $R^{op}$ is a ring such that $R^{op}=R$ (as additive groups) but the multiplication in $R^{op}$ is given by the rule $a\cdot b = ba$.
 Recall that $\pCCR$ and $\CC_R'$ are the sets of left and right regular elements of the ring $R$, respectively, and $S_l(R)$ and $S_r(R)$ are the largest left and right Ore sets of $R$ that consists of regular elements of $R$.

\begin{lemma}\label{a8Mar15}
Let $R$ be a ring.
\begin{enumerate}
\item In the set $\pCCR$ there exists the largest (w.r.t. inclusion) left denominator set of $R$, denoted by $\pSlR$. The set $S_r(R)$ is the largest (w.r.t. inclusion) right denominator set of $R$ in $\pCCR$.
\item In the set $\CC_R'$ there exists the largest (w.r.t. inclusion) right denominator set of $R$, denoted by $S_r'(R)$. The set $S_l(R)$ is the largest (w.r.t. inclusion) left denominator set of $R$ in $\CC_R'$.
\end{enumerate}
\end{lemma}

{\it Proof}. 1. If $S$ and $T$ are left denominator sets of the ring $R$ such that $S,T\subseteq \pCCR$. The submonoid, denoted by $ST$,  of $\pCCR$ that is  generated by $S$ and $T$ does not contain $0$. By \cite[Lemma 2.4.(2)]{Bav-stronglquot}, $ST$ is a left denominator set of $R$. Hence, the set $\pSlR$ exists and is the union of all left denominator sets of $R$ in $\pCCR$.

 If $D$ is a right denominator set of $R$ in $\pCCR$ then $ \ass_R(D)=0$, and so $D\subseteq \CC_R$. Therefore, $S_r(R)$ is the largest right denominator set of $R$ in $\pCCR$.

 2. Statement 2 follows from statement 1 (by applying statement 1 to the opposite ring $R^{op}$ of $R$).
 $\Box $

 {\it Definition}. The set $\pSlR$ is called the {\em largest left regular left denominator set} of $R$ and the ring $\pQlR : = \pSlR^{-1}R$ is called the {\em left regular left quotient ring} of $R$. Similarly, the set $S_r'(R)$ is called the {\em largest right regular right denominator set} of $R$ and the ring $Q_r'(R):= RS_r'(R)^{-1}$ is called the {\em right regular right quotient ring} of $R$.

If $\pSlR= \pCCR$ then $\pQlR = \pQ_{l,cl}(R)$. If $S_r'(R)= \CC_R'$ then $Q_r'(R) = Q_{r,cl}'(R)$.

The next lemma shows that if the ring $Q_l(R)$ is a left Artinian ring (respectively, semisimple Artinian) ring then so is the ring $\pQ_l(R)$. The reverse implication is usually wrong. For example, in the case of the algebra $\mI_1= K \langle x, \frac{d}{dx}, \int\rangle$ of the polynomial integro-differential operators over a field $K$ of characteristic zero, the ring $Q_l(\mI_1)$ is neither left nor right Noetherian ring and not a domain (see \cite{larglquot}) but the ring $\pQ_l(\mI_1)$ is a division ring and
 $\pQ_l(\mI_1)=\pQ_{l,cl}(\mI_1)$ (Theorem \ref{19Mar15}.(1)).

\begin{lemma}\label{b8Mar15}
Let $R$ be a ring.
\begin{enumerate}
\item If the ring $Q_l(R)$ is a left Artinian ring then $S_l(R) = \CC_R = \pCCR = \pSlR$ and $Q_l(R) = Q_{l,cl}(R) = \pQ_{l,cl}(R) = \pQlR$ is a left Artinian ring.
\item If the ring $Q_l(R)$ is a semisimple  Artinian ring then $S_l(R) = \CC_R = \pCCR = \pSlR$ and $Q_l(R) = Q_{l,cl}(R) = \pQ_{l,cl}(R) = \pQlR$ is a semisimple Artinian ring.
\end{enumerate}
\end{lemma}

{\it Proof}. 1.If $Q_l(R)$ is a left Artinian ring then $\pCCR\subseteq Q_l(R)^*$. By \cite[Theorem 2.8.(1)]{larglquot}, $S_l(R)= R\cap Q_l(R)^*$. By intersecting with $R$ the following inclusions of subsets of  the ring $Q_l(R)$, $S_l(R)\subseteq \CC_R\subseteq \pCCR \subseteq Q_l(R)^*$ and $S_l(R) \subseteq \pSlR \subseteq \pCCR \subseteq Q_l(R)^*$, we obtain the inclusions  $S_l(R)\subseteq \CC_R\subseteq \pCCR \subseteq S_l(R)$ and $S_l(R) \subseteq \pSlR \subseteq \pCCR \subseteq S_l(R)$. Therefore, $S_l(R) = \CC_R = \pCCR = \pSlR$ and $Q_l(R) = Q_{l,cl}(R) = \pQ_{l,cl}(R) = \pQlR$ is a semisimple Artinian ring.

2. Statement 2 follows from statement 1.
 $\Box $

{\bf Semisimplicity criteria for the ring $\pQ_l(R)$}. In general, the question of existence of the ring $\pQ_{l,cl}(R)$ is a difficult one. In general, the ring $\pQ_{l,cl}(R)$ does not exits but the ring $\pQ_l(R)$ always does. If the ring $\pQ_{l,cl}(R)$ exists  then $\pQ_{l,cl}(R)=\pQ_l(R)$. The next theorem states that if the ring $\pQ_l(R)$ is a left Artinian ring or a semisimple Artinian ring then so is the ring $\pQ_{l,cl}(R)$, and vice versa.

\begin{theorem}\label{8Mar15}
Let $R$ be a ring. Then
\begin{enumerate}
\item $\pQlR$ is a left Artinian ring iff $\pQ_{l,cl}(R)$ is a left Artinian ring. If one of the equivalent conditions holds then $\pSlR = \pCCR$ and $\pQlR = \pQ_{l, cl}(R)$.
\item $\pQlR$ is a semisimple  Artinian ring iff $\pQ_{l,cl}(R)$ is a semisimple Artinian ring. If one of the equivalent conditions holds then $\pSlR = \pCCR$ and $\pQlR = \pQ_{l, cl}(R)$.

\end{enumerate}
\end{theorem}

{\it Proof}. 1. $(\Rightarrow )$  Let  $\pS := \pSlR$, $\pga := \ass_R(\pS )$ and ${}'\pi : R\ra R/ \pga$, $r\mapsto \br := r+\pga$.

(i) ${}'\pi (\pCCR ) \subseteq \pCC_{R/\pga }$: Suppose that $\br \bc =0$ for some elements $r\in R$ and $c\in \pCCR$. We have to show that $r\in \pga$. The element $a:= rc$ belongs to the ideal $\pga$. Then $0=sa=src$ for some element $s\in \pS$, and so $sr=0$ (since $c\in \pCCR $). Therefore, $r\in \pga$.

(ii) $\pCC_{R/ \pga }= \CC_{R/ \pga }$: Since $\pS \in \Den_l(R, \pga )$, we have the inclusion ${}'\pi (\pS ) \in \Den_l(R/ \pga , 0)$ and the ring $\pQlR \simeq {}'\pi (\pS )^{-1} (R/ \pga )$ is a left  Artinian ring. Hence, $\pCC_{R/ \pga } = \CC_{R/\pga }$.

(iii) $\ga := \ass_R(\pCCR ) = \pga$: The inclusion $\pS\subseteq \pCCR$ implies the inclusion $\pga \subseteq \ga$. Then, by the statements (i) and (ii), ${}'\pi(\pCCR ) \subseteq \CC_{R/ \pga}$, and so  we must have $\pga = \ga$.

(iv) $\pCCR \in \Ore_l(R)$: Let $c\in \pCCR$ and $r\in R$. By the statements (i) and (ii), ${}'\pi (\pCCR ) \subseteq \CC_{R/ \pga }$. Hence, $\bc\in \CC_{R/ \pga }$. By the assumption, the ring $\pQlR \simeq {}'\pi (\pS)^{-1} ( R/ \pga )$ is a left Artinian ring. Hence, $Q_{l, cl} (R/ \pga )\simeq \pQlR$ is a left Artinian ring and the elements of the set $\CC_{R/ \pga }$ are units in the ring $\pQlR$. In particular, the element $\bc$ is so. So, $\br\bc^{-1} = \bs^{-1} \oa$ for some elements $s\in \pS$ and $a\in R$. Then $sr-ac\in \pga$, and so $s'(sr-ac)=0$ for some elements $s'\in \pS$. So, $s's\cdot r= s'a\cdot c$ where $s's\in \pS \subseteq \pCCR$, and the statement (iv) follows.

(v) $\pCCR \in \Den_l(R, \pga )$: This follows from the statements (iii) and (iv).

(vi) $\pS = \pCCR$ (by the maximality of $\pS$) and so $\pQlR = \pQ_{l, cl}(R)$ is a left Artinian  ring.

$(\Leftarrow )$ This implication is obvious.

2. Statement 2 follows from statement 1.  $\Box $

In view of Theorem \ref{8Mar15}, {\em all the semisimplicity criteria for the ring $\pQ_{l,cl}(R)$ are also semisimplicity criteria for the ring $\pQ_l(R)$, and vice versa}.

{\bf The canonical homomorphism $\phi : Q_l(R)\ra \pQ_l(R)$}. The next proposition shows that there is a canonical ring homomorphism $\phi : Q_l(R)\ra \pQ_l(R)$ and gives a criterion for $\phi$ to be an isomorphism.
\begin{proposition}\label{1Mar15}
Let $R$ be a ring. Then
\begin{enumerate}
\item $S_l(R)\subseteq \pSlR \subseteq \pCCR$, and so $\ass_R (S_l(R))\subseteq \ass_R (\pSlR )\subseteq \ass_R (\pCCR )$.
\item The map $\phi : Q_l(R)\ra \pQ_l(R)$, $s^{-1}r\mapsto s^{-1}r$, is a ring $R$-homomorphism with kernel $S_l(R)^{-1}\pga$ where $\pga = \ass_R(\pSlR )$.
    \item $\phi$ is an isomorphism iff $\pga =0$ iff $S_l(R)= \pSlR$.
    \item The rings $Q_l(R)$ and $\pQ_l(R)$ are $R$-isomorphic iff one of the equivalent conditions of statement 3 holds.
\end{enumerate}
\end{proposition}

{\it Proof}. 1. Notice that  $S_l(R)\subseteq \CC_R\subseteq \pCCR$. By Lemma \ref{a8Mar15}.(1), $S_l(R) \subseteq \pSlR \subseteq \pCCR$.

2. Statement 2 follows from statement 1: By statement 1, the map $\phi$ is well-defined.  If $\phi (s^{-1}r) =0$ then $\frac{r}{1}=0$ in $\pQ_l(R)$, and so $r\in \pga$, i.e. $\ker (\phi ) = S_l(R)^{-1}\pga$.

3. $\phi$ is an isomorphism iff $\ker (\phi ) =0$ and $\phi$ is a surjection iff $\pga =0$ (since $S_l(R)\subseteq \CC_R$) and $\phi$ is a surjection iff $S_l(R)= \pSlR$ and $\phi$ is a surjection iff $S_l(R)= \pSlR$ iff $\pga =0$.

4. We have to show that $(3\Leftrightarrow 4)$. The implication $(3\Rightarrow 4)$
 is obvious. Conversely, suppose that $\v : Q_l(R)\ra \pQ_l(R)$ is an $R$-isomorphism ($\v (rq) = r\v (q)$ for all $r\in R$ and $q\in Q_l(R)$).  Then $R\subseteq Q_l(R)$ and $\pga \subseteq \ker (\v ) =0$, i.e. $\pga =0$. So, statement 3 holds.   $\Box $

 The next corollary shows that if $S_l(R)\neq \pSlR$ then the ring $Q_l(R)$ is not  left Artinian.

 \begin{corollary}\label{a9Mar15}

\begin{enumerate}
\item Let $R$ be a ring such that $S_l(R)\neq \pSlR$ or, equivalently, $\pga := \ass_R(\pSlR )\neq 0$ (Proposition \ref{1Mar15}.(3)). Then $Q_l(R)$ is not a left Artinian  ring. In particular, $Q_l(R)$ is not a semisimple Artinian ring.
\item If, in addition, the ring $R$ is a $K$-algebra over a field $K$ then, for all algebras $A$, $Q_l(R\t A)$ is not a left Artinian ring. In particular,  $Q_l(R\t A)$ is not a semisimple Artinian ring.
\end{enumerate}
\end{corollary}

{\it Proof}. 1. Suppose that the ring $Q_l(R)$ is a left Artinian ring. Then, by Lemma \ref{b8Mar15}.(1), $\pSlR = \CC_R$, and so $\pga =0$, a contradiction.

2. Clearly, $\pCCR \subseteq \pCC_{R\t A}$. Hence, $\pSlR\subseteq \pS_{R\t A}$ (by Lemma \ref{a8Mar15}.(1)), and so $0\neq \pga = \ass_R(\pS_l)\subseteq \ass_{R\t A} (\pS_{R\t A})$. By statement 1, $Q_l(R\t A)$ is not a left Artinian ring.
 $\Box $

{\bf Applications to the algebras of polynomial integro-differential operators}.
 Let $K$ be a field of characteristic zero, $K[x]$ be a polynomial algebra in a single variable $x$, $\der := \frac{d}{dx}$ and $\int : K[x]\ra K[x]$, $x^n\mapsto \frac{x^{n+1}}{n+1}$ for all $n\geq 0$ be the integration. The following subalgebras of $\End_K(K[x])$,  $A_1= K\langle x,\der \rangle$ and $\mI_1=K\langle x,\der , \int\rangle$,  are called the {\em first Weyl algebra} and the {\em algebra of polynomial integro-differential operators}, respectively. By definition, $A_n:= A_1^{\t n}$ is called  the $n$'th {\em Weyl algebra} and $\mI_n := \mI_1^{\t n}$ is called the {\em algebra of polynomial integro-differential operators}. The Weyl algebra $A_n$ is a Noetherian domain, and so $Q_{l,cl}(A_n)$ is a division ring. For the algebra $\mI_n$, neither the ring $Q_{l,cl}(\mI_n)$ nor the ring $Q_{r,cl}(\mI_n)$  exists (Lemma \ref{b21Mar15}.(1)). The rings $Q_l(\mI_n)$ and $Q_r(\mI_n)$ are neither left nor right Noetherian and not domains (Lemma \ref{b21Mar15}.(2)).

As an easy  application of Corollary \ref{a9Mar15} we have  the next result. A more strong result of that kind is Lemma \ref{b21Mar15} where different arguments are used in its proof.
\begin{corollary}\label{b9Mar15}
For all $n \geq 1$, the rings $Q_l(\mI_n)$ (resp., $Q_r(\mI_n)$) are not left (resp., right) Artinian. Moreover, for all algebras $A$, the rings $Q_l(\mI_n\t A)$ (resp., $Q_r(\mI_n\t A)$)  are not left (resp., right)  Artinian.
\end{corollary}

{\it Proof}. The ring $\mI_n$ is isomorphic to its opposite ring \cite{larglquot}. In view of this fact and Corollary \ref{a9Mar15}, it suffices to show that $S_l(\mI_1) \neq \pS_l(\mI_1)$. The set $S_\der :=\{ \der^i\, | \, i\geq 0\}$ is a left denominator set of the algebra $\mI_1$ (see \cite{larglquot}) with $\ass_{\mI_1}(S_\der)\neq 0$ since $\der\cdot (1-\int\der ) = \der - \der\int\der = \der -1\cdot \der =0$. Clearly, $\der\in \pCC_{\mI_n}$ since $\der\int =1$. Then, $\der\in \pS_l(\mI_1)\backslash S_l(\mI_1)$, as required.  $\Box $

By Proposition \ref{1Mar15}, the ring homomorphism $\phi$ is the composition of the following ring homomorphisms:
\begin{equation}\label{QQl1}
\phi : Q_l(R)\stackrel{\pi'}{\ra} \bQ := Q_l(R)/S_l(R)^{-1}\pga \ra \pQ_l(R)\simeq T^{-1}\bQ
\end{equation}
where $\pi' (a) = a+S_l(R)^{-1}\pga$ and $T\in \Den_l(\bQ )$ is the multiplicative subset of $(\bQ , \cdot )$ generated by the group of units $\bQ^*$ of the ring $\bQ$ and the set $\pi' (\pSlR )$.

\begin{corollary}\label{a1Mar15}
Let $R$ be a ring. Suppose that $\CP$ is a property of rings that is preserved by left localizations and  passing to factor ring. If the ring $Q_l(R)$ satisfies the property $\CP$ then so does the ring $\pQ_l(R)$. In particular, if the ring $Q_l(R)$ is a semisimple (respectively, left Artinian; left Noetherian) then so is the ring $\pQ_l(R)$.
\end{corollary}

{\it Proof}. The corollary follows from (\ref{QQl1}). $\Box $






The next lemma gives plenty of examples of algebras for which neither left nor right classical quotient ring exists. This is true for the algebras $\mI_n$.
\begin{lemma}\label{b21Mar15}
Let $A$ be an algebra over $K$.
\begin{enumerate}
\item The set $\CC_{\mI_1\t A}$ of regular elements of the algebra $\mI_n\t A$ is neither a left nor right Ore set. Therefore, the rings $Q_{l,cl}(\mI_n\t A)$ and $Q_{r,cl}(\mI_n\t A)$ do not exist.
\item The algebras $Q_l(\mI_n\t A)$ contain  infinite direct sums of nonzero left ideals and so they are  not  left Noetherian algebras. Similarly, the algebras $Q_r(\mI_n\t A)$ contain  infinite direct sums of nonzero right ideals and so they are not right  Noetherian algebras.
\end{enumerate}
\end{lemma}

{\it Proof}. 1.  Clearly, $\mI_n\t A= \mI_1\t\mI_{n-1}\t A$ and $\CC_{\mI_1}\subseteq \CC_{\mI_n\t A}$.  Let $e_{00}:= 1-\int\der$. The element $a:=\der_1+\int_1$ belongs to the set $\CC_{\mI_1}$, $\mI_1e_{00}\cap \mI_1a=0$ and $e_{00}\mI_1\cap a\mI_1=0$, see the proof of \cite[Theorem 9.7]{Bav-intdifline}. Hence, $a\in \CC_{\mI_n\t A}$ and
$(\mI_n\t A) e_{00}\cap (\mI_n\t A)a = (\mI_1\t \mI_{n-1}\t A) e_{00} \cap (\mI_1\t \mI_{n-1}\t A) a= (\mI_1e_{00}\cap \mI_1a)\t \mI_{n-1}\t A=0$, and similarly $e_{00}(\mI_n\t A) \cap a
(\mI_n\t A)=0$. This means that $\CC_{\mI_n\t A} \notin \Ore_l( \mI_n\t A) \cup \Ore_r( \mI_n\t A)$.

2. The algebra $\mI_n$ contains infinite direct sums of nonzero left ideals \cite{Bav-intdifline}, hence so do the algebras $\mI_n\t A$, and statement 2 follows.  $\Box $


\section{Properties of $\pSlR$ and $\pQ_{l,cl}(R)$}\label{S5MAR}

In this section, some properties of $\pS_l(R)$ and $\pQ_{l,cl}(R)$ are established (Theorem \ref{9Mar15}). The main motive is to develop practical tools for finding the ring $\pQ_{l,cl}(R)$. A key idea is that in order to find $\pQ_{l,cl}(R)$ there is no need to know explicitly the set $\pS_l(R)$. It suffices to replace $\pS_l(R)$ with another left denominator set that yields the same result, see Theorem \ref{9Mar15}.(5). Further developing this idea sufficient conditions are found for the ring $\pQ_{l,cl}(R)$ to be isomorphic to $Q_{l,cl}(R/ \ass_R(\pS_l(R)))$ (Theorem \ref{11Mar15}).

\begin{lemma}\label{a5Mar15}
Let $R$ be a ring, $S$ be a multiplicative subset of $\pCCR$ such that $\ga' := \ass_R(S)$ is an ideal, $\pi': R\ra \bR':= R/\ga'$, $r\ra r+\ga'$. If $\pi'(S) \in \Ore_l(\bR')$ then $S\in \Den_l(R, \ga')$.
\end{lemma}

{\it Proof}. Since $S\subseteq \pCCR$, it suffices to show that $S\in \Ore_l(R)$. Given $s\in S$ and $r\in R$. Then $\bs_1\br = \br_1\bs$ for some elements $ s_1\in S$ and $r_1\in R$ (since $\pi'(S)\in \Ore_l(\bR')$). Then $s_1r-r_1s\in \ga'$ and so $s_2(s_1r-r_1s)=0$ for some elements $s_2\in S$, and we are done (since $ s_2s_1\cdot r = s_2r\cdot s$).
 $\Box $

\begin{proposition}\label{b5Mar15}
Let $R$ be a ring, $S$ be a multiplicative subset of $\pCCR$ such that $\ga' := \ass_R(S)$ is an ideal of $R$, $\pi': R\ra \bR':= R/\ga'$, $r\ra r+\ga'$. We keep the notation of Theorem \ref{3Mar15}. Then
\begin{enumerate}
\item $S\subseteq \{ c\in R\, | \, \cdot c_{\ga'}$ and $\cdot c_{R/ \ga'}$ are injections$\}$. If, in addition,  $S= \pCCR$ then $\pCCR = \{ c\in R\, | \,  \cdot c_{\ga}$ and $\cdot c_{R/ \ga}$ are injections$\}$ (where $\ga =\ass_R(\pCCR )$ is an ideal of $R$, by the assumption).
\item $\pi'(S)\subseteq \CC_{\bR'}$. In particular, $\pi (\pCCR ) \subseteq \CC_{\bR}$ provided $\ga$ is an ideal of $R$.
\item $S\subseteq \{ c\in R\, | \, \pi'(c)\in \CC_{\bR'}$,  $\cdot c_{\ga'}$ is an  injection$\}$.  If $\ga$ is an ideal of $R$ then  $\pCCR = \{ c\in R\, | \,  \pi (c)\in \CC_{\bR}$,  $\cdot c_{\ga}$ is an  injection$\}$.
\item  $S\in \Den_l(R, \ga')$  iff $\pi'(S)\in \Den_l(\bR', 0)$ iff $\pi'(S)\in \Ore_l(\bR')$ iff $S\in \Ore_l(R)$.
 \item If $\ga$ is an ideal of $R$ then  $\pCCR\in \Den_l(R, \ga )$  iff $\pi (\pCCR )\in \Den_l(\bR , 0)$ iff $\pi (\pCCR )\in \Ore_l(\bR)$ iff $\pCCR \in \Ore_l(R)$.
\end{enumerate}
\end{proposition}

{\it Proof}. 1. Clearly, $T:= \{ c\in R\, | \, \cdot c_{\ga'}$ and $\cdot c_{R/\ga'}$ are injections$\}\subseteq \pCCR$. Given $c\in S$. In order to show that $S\subseteq T$,  we have to prove that $s\in T$. Since $c\in S\subseteq \pCCR$, the map $ \cdot c_{\ga'}$ is an injection. It remains to show that $ \cdot c_{R/ \ga'}$ is also an injection. If $ \br \bc =0$ then $rc\in \ga'$, and so $src=0$ for some $s\in S$. Then $sr=0$ (since $c\in S\subseteq \pCCR$), and so $r\in \ga'$, i.e. $\br =0$. Therefore, $S\subseteq T$.

If $S= \pCCR$ then $\pCCR = T$ (since $T\subseteq \pCCR$).

2. By statement 1, $\pi'(S)\subseteq \pCC_{\bR'}$. To prove that the inclusion $\pi'(S)\subseteq \CC_{\bR'}$ holds it remains to show that each element $\bs : = s+\ga'$ (where $s\in S$) is a right regular element of the ring $\bR'$. Suppose that $\bs \br =0$ for some element $r\in R$. Then $sr\in \ga'$, and so $tsr=0$ for some element $t\in S$. This implies that $ r\in \ga'$, i.e. $\br =0$, and so $\bs$ is a right regular element of $R$.

3. Statement 3 follows from statements 1 and 2.

4. Since $ S\subseteq \pCCR$ and $\ga' = \ass_R(S)$, $S\in \Den_l(R, \ga')$ iff $S\in \Ore_l(R)$. By statement 2, $\pi'(S)\subseteq \CC_{\bR'}$. So, $\pi'(S)\in \Den_l( \bR', 0)$ iff $\pi'(S)\in \Ore_l(\bR )$. It remains to show that the first `iff' holds.

Suppose that $S\in \Den_l(R, \ga')$. Then $\pi'(S) \in \Den_l(\bR , 0)$.

Suppose that $\pi'(S)\in \Den_l(\bR , 0)$. By Lemma \ref{a5Mar15}, $S\in \Den_l(R, \ga')$.

5. Statement 5 is a particular case of statement 4 where $S= \pCCR$.  $\Box $

For a ring $R$, let $\pAsslR := \{ \ass_R(S)\, | \, S\in \Den_l(R), S\subseteq \pCCR \}$. The set $(\pAsslR , \subseteq )$ is a poset.

\begin{theorem}\label{9Mar15}
Let $R$ be a ring and $\pga :=\ass_R(\pSlR )$. Then
\begin{enumerate}
\item $\pSlR = \pCCR\cap (\pSlR +\pga )$.
\item $\pga$ is the largest element in $\pAsslR$.
\item $\pSlR$ is a maximal element (w.r.t. inclusion) in the set $\{ S\in \Den_l(R)\, | \, S= \pCCR \cap (S+\ass_R(S))\}$.
    \item $\pSlR +\pga \in \Den_l(R, \pga )$.
    \item $\pQ_l(R)\simeq (\pSlR +\pga )^{-1}R$.
\end{enumerate}
\end{theorem}

{\it Proof}. 1. Let $\pS = \pSlR$ and $T:= \pCCR\cap (\pSlR +\pga )$.  Clearly, $\pS \subseteq T$. In order to prove that $\pS \supseteq T$, it suffices to show that $T\in \Den_l(R)$ (since $T\subseteq \pCCR$ and  by maximality of $\pS$ (Lemma \ref{a8Mar15}.(1)), $\pS \supseteq T$).  By the very definition, $T$ is a multiplicative set in $R$ such that ${}'\pi(T) = {}'\pi (\pS)$ where ${}'\pi : R\ra R/ \pga $, $r\mapsto \br := r+\pga$. Since $T\subseteq \pCCR$, $T\in \Den_l(R)$ iff $T\in \Ore_l(R)$. Let us show that  $T\in \Ore_l(R)$. Given elements $t\in T$ and $r\in R$. Then $\bt \in {}'\pi (T) = {}'\pi (\pS )$. Since ${}'\pi (\pS ) \in \Den_l(R/ \pga , 0)$, $\bs \br = \br_1\bt$ for some elements $s\in \pS$ and $r_1\in R$. Then $sr-r_1t\in \pga$, and so $s'(sr-r_1t)=0$
 for some element $s'\in \pS$, i.e. $s's\cdot r=s'r\cdot t$ where $s's\in \pS\subseteq T$. Therefore, $T\in \Ore_l(R)$.

2. Statement 2 follows from the maximality of $\pS$ (Lemma \ref{a8Mar15}.(1)).

3.  Every element $S\in \Den_l(R)$ such that $S= \pCCR \cap (S+\ass_R(S))$ consists of left regular elements, hence $S\subseteq \pS$. Now, statement 3 follows from statement 1.

4.  Let $T':= \pSlR +\pga$. Then ${}'\pi (T') = {}'\pi (\pS )$, and so the set $T'$ is a multiplicative set of $R$.

(i) $T'\in \Ore_l(R)$:  Given elements $t\in T'$ and $r\in R$. Then $\bt \in {}'\pi (T') = {}'\pi (\pS )$. Since ${}'\pi (\pS ) \in \Den_l(R/ \pga , 0)$, $\bs \br = \br_1\bt$ for some elements $s\in \pS$ and $r_1\in R$. Then $sr-r_1t\in \pga$, and so $s'(sr-r_1t)=0$
 for some element $s'\in \pS$, i.e. $s's\cdot r=s'r\cdot t$ where $s's\in \pS\subseteq T'$. Therefore, $T'\in \Ore_l(R)$.

(ii) $\ass_R(T') = \pga$: The inclusion $\pS\subseteq T'$ implies the inclusion $\pga \subseteq \gb :=\ass_R(T')$. If $r\in \gb$, i.e. $tr=0$ for some element $t=s+a\in T'$ where $s\in \pS$ and $a\in \pga$. Fix an element $s'\in \pS$ such that $s'a=0$. Then $0=s'tr= s'(s+a)r= s'sr$, and so $r\in \pga$ since $s's\in \pga$.

(iii)  $T'\in \Den_l(R)$: We have to show that if $rt=0$ for some elements $r\in R$ and $t\in T'$ then $r\in \pga$. The element $t\in T'$ is a sum $s+a$  where $s\in \pS$ and $a\in \pga$. Then the equality $0=rt= r(s+a)$ can be written as $rs= -ra\in \pga$. Hence, $s'rs=0$ for some element $s'\in \pS$, and so $s' r=0$ (since $s\in \pS\subseteq \pCCR$).  Therefore, $r\in \pga$ (as $s'\in \pS$).

5. By statement 4, $\pQ_l(R)\simeq {}'\pi(\pS )^{-1} (R/ \pga )= {}'\pi (\pS+\pga )^{-1}(R/ \pga )\simeq (\pS+\pga )^{-1}R$.  $\Box $

{\bf Sufficient conditions for $\pQ_{l, cl}(R)\simeq Q_{l,cl}(R/\ass_R(\pCCR ))$}.
 The next theorem gives sufficient conditions for the ring $\pQ_{l, cl}(R)$ to be isomorphic to the ring $Q_{l,cl}(R/\ass_R(\pCCR ))$.

\begin{theorem}\label{11Mar15}
Let $R$ be a ring, $\pCC = \pCCR$ and $\ga = \ass_R(\pCC )$. Suppose that $\ga$ is an ideal of the ring $R$, the set $\pbCC := \pi (\pCC )$ is a dense subset of $\CC_{\bR }$ in $\bR := R/\ga$ where $\pi : R\ra \bR $, $r\mapsto \br =r+\ga$, and $\CC_{\bR }\in \Ore_l(R)$.  Then
\begin{enumerate}
\item $\pCC \in \Den_l(R, \ga )$ and $\pbCC\in \Den_l(\bR , 0)$.
\item $\pi^{-1} (\CC_{\bR})\in \Den_l(R, \ga )$ and $\pCC \subseteq \pi^{-1} (\CC_{\bR})$.
\item $\pQ_{l,cl} (R)\simeq \pbCC^{-1} \bR \simeq \pi^{-1} (\CC_{\bR})^{-1} R\simeq Q_{l, cl}(\bR )$.
\end{enumerate}
\end{theorem}

{\it Remark}. If $\ga$ is an ideal of the ring $R$ then $\pbCC\subseteq \CC_{\bR}$, by Proposition \ref{b5Mar15}.(2).

{\it Proof}. 1. (i) $\pCC \in \Ore_l(R)$: Clearly, $\pCC$ is a multiplicative set of $R$. Given elements $c\in \pCC$ and $r\in R$. We have to find elements $c'\in \pCC$ and $r'\in R$ such that $c' r= r'c$. By the assumption, $\CC_{\bR}\in \Ore_l(\bR )$. So, let $\bQ := Q_{l, cl}(\bR )$. Since $\pbCC\subseteq \CC_{\bR}$, we have $\br \bc^{-1} = \bs^{-1} \br_1$ for some elements $\bs = s+\ga\in \CC_{\bR }$ (where $s\in R$) and $r_1\in R$. We can write $\bs_1\br = \br_1\bc$. By the assumption,  the set $\pbCC$ is a dense subset of $\CC_{\bR}$. Fix an element $\br_2\in \bR$ (where $r_2\in R$) such that $\bc_1:= \br_2\bs \in \pbCC$ for some element $c_1\in \pCC$. Then the equality $\br_2\bs_1\br = \br_2\br_1\bc$ can be written as $\bc_1\br = \br_2\br_1\bc$. Hence, $c_1r-r_2r_1c\in \ga$, and so there exists an element $c_2\in \pCC$ such that $c_2(c_1r-r_2r_1c)=0$. Notice that  $c':= c_2c_1\in \pCC$, $r' := c_2r_2r_1\in R$ and $ c'r = r'c$.

(ii) $\pCC \in \Den_l(R, \ga )$: The inclusion follows from Lemma \ref{a5Mar15} and (i).

 By (ii), $\pbCC\in \Den_l(\bR , 0)$ and $\pQ_{l,cl} (R) \simeq \pbCC^{-1} \bR$.

2. By statement 1, $\pCC\in \Den_l(R, \ga )$ and $\pCC\subseteq \pi^{-1} (\pCC_{\bR})$. By the assumption, $\CC_{\bR}=\pi(\pi^{-1}(\CC_{\bR}))\in \Den_l(\bR , 0)$. Therefore, $\pi^{-1} (\CC_{\bR})\in \Den_l(R, \ga )$.

3. By statement 2, $\pi^{-1}(\CC_{\bR })^{-1} R\simeq \CC_{\bR}^{-1} \bR = \bQ$. By the assumption, the set $\pbCC$ is dense in $\CC_{\bR}$. By statement 1,  $\pbCC\in \Den_l(\bR , 0)$.  Hence, $\pQ\simeq \pbCC^{-1} \bR \simeq \pQ_{l,cl}(R)$.  $\Box $


\section{The classical left regular left quotient ring of the algebra of polynomial integro-differential operators $\mI_1$}\label{CLRLQIN}

The aim of this section is to find for the  algebra of polynomial integro-differential operators $\mI_1$ its classical left regular left quotient ring $\pQ :=\pQ_{l,cl}(\mI_1)$
 and classical right  regular right quotient ring $Q':=\pQ_{r,cl}(\mI_1)$, and to show that both of them are canonically isomorphic to to the classical quotient ring of the  Weyl algebra $A_1$ (Theorem \ref{19Mar15}). The sets $\pCC_{\mI_1}$ and $\CC_{\mI_1}'$ are described in Theorem \ref{24Mar15}.  The algebra $A_1$ is a Noetherian domain so $Q(A_1):=Q_{l,cl}(A_1)\simeq Q_{r,cl}(A_1)$. The key idea in finding the rings $\pQ$ and $Q'$ is to use Theorem \ref{28Feb15}. The most difficult part is to verify that the set $\pbCC$ is dense in $\CC_{\mI_n}$ (see, Corollary \ref{b19Mar15}).

We start this section with collecting necessary facts about the algebra $\mI_1$ that are used in the proofs (their proofs are given in \cite{Bav-intdifline}).

{\bf The algebra $\mI_1$ of polynomial integro-differential operators}.  Let us recall some of the properties of the algebra  $\mI_1$.  Let $K$ be a field of characteristic zero, $P_1=K[x]$  and $E_1=\End_K(P_1)$ be the algebra of all $K$-linear maps from $P_1$ to $P_1$. Recall that the algebra $\mI_1$ of polynomial integro-differential operators is the subalgebra of $E_1$ generated by the elements $x$, $\der = \frac{d}{dx}$ and $\int $. The algebra $\mI_1$ contains the Weyl algebra $A_1=K\langle x,\der\rangle$. The algebra $A_1$ is a Noetherian domain but  the algebra $\mI_1$ is neither left nor right Noetherian domain. Moreover, it contains infinite direct sums of nonzero left and right ideals. Neither left nor right classical quotient ring exits (Lemma \ref{b21Mar15}.(1)). The largest left quotient ring $Q_l(\mI_1)$ and the largest right quotient ring $Q_r(\mI_1)$  are neither left nor right Noetherian rings (Lemma \ref{b21Mar15}.(2)). The algebra $\mI_1$ admits a single proper ideal $F=\oplus_{i,j\in \N} Ke_{ij}$ where $e_{ij} = \int^i \der^j-\int^{i+1}\der^{j+1}$, and ${}_{K[\int]}F_{K[\der]}=K[\int]e_{00}K[\der]\simeq K[\int]\t K[\der]$.
 The factor algebra $\mI_1/ F$ is isomorphic to the algebra $A_{1,\der}$ which is a the localization of the Weyl algebra $A_1$ at the powers of the element $\der$. Each element $a\in \mI_1$ is a unique sum
\begin{equation}\label{acan}
a= \sum_{i>0}a_{-i}\der^i + a_0+\sum_{i>0}\int^ia_i+\sum_{i,j} \l_{ij} e_{ij}
\end{equation}
where $a_k\in K[H]$, $H:= \der x$ and $\l_{ij}\in K$.  Since $\der \int =1$, we have the equalities $\der e_{ij} = e_{i-1, j}$, $e_{ij}\der = e_{i,j+1}$, $\int e_{ij} = e_{i+1, j}$ and $e_{ij}\int = e_{i,j-1}$ (where $e_{-1,j}:=0$ and $e_{i,-1}:=0$).
 The algebra $\mI_1$  is generated by the elements $\der $, $H:=
\der x$ and $\int$ (since $x=\int H$) that satisfy the defining
relations: $$\der \int = 1,
\;\; [H, \int ] = \int, \;\; [H, \der ] =-\der , \;\; H(1-\int\der
) =(1-\int\der ) H = 1-\int\der .$$
Since the algebra $\mI_1/F$ is a domain, $\pCC_{\mI_1} = \{ a\in \mI_1\, | \, \ker (\cdot a_F)=0\}$ where $\cdot a_F: F\ra F$, $f\mapsto fa$. The right $\mI_1$-module $F$ is the direct sum $\oplus_{i\in \N}e_{i0}K[\der ]$ of isomorphic right $\mI_1$-modules. The right $\mI_1$-module $e_{00}K[\der ]$ is a free right $K[\der]$-module of rank 1. When we identify the right $K[\der ]$-modules $e_{00}K[\der]$ and $P_1':=K[\der]$, the right $\mI_1$-module structure on the polynomial algebra $P_1'$ is given by the rule: For $i\geq 0$, $\der^i \cdot \der = \der^{i+1}$, $\der^i\cdot \int = \der^{i-1}$ $(i\geq 1)$ and $1\cdot \int =0$, $\der^i \cdot H = \der^i (i+1)$. So, $\pCC_{\mI_1} = \{ a\in \mI_1 \, | \, \ker ( \cdot a_{P_1'})=0\}$.

The algebra
$\mI_1$ admits the involution $*$ over the field $K$:
$\der^* =\int $, $\int^* = \der$  and $H^* = H,$ i.e.  it is a $K$-algebra
{\em anti-isomorphism} ($(ab)^* = b^* a^*$) such that $a^{**} =a$.
Therefore, the algebra $\mI_1$ is {\em self-dual}, i.e.  it is
isomorphic to its opposite algebra $\mI_1^{op}$. As a result, the
left and right properties of the algebra $\mI_1$ are the same.
Clearly, $e_{ij}^* = e_{ji}$ for all $i,j\in \N$, and so $F^* =
F$.

\begin{lemma}\label{c16Mar15}
  Suppose that $T\in \Den_l(R)$ and $S$ be a multiplicative set of $R$ such that $S\subseteq T$, $\ass_R(S) = \ass_R(T)$ and for each element $t\in T$ there exists an element $r\in R$ such that $rt\in S+\ass_R(T)$. Then $S\in \Den_l(R)$ and $S^{-1}R\simeq T^{-1}R$.
\end{lemma}

{\it Proof}. (i) $S\in \Ore_l(R)$: Given elements $s\in S$ and $r\in R$. We have to show that $s' r = r' s$ for some elements $s'\in S$ and $r'\in R$. Since  $T\in \Ore_l(R)$ and $s\in S\subseteq T$, $tr=r_1s$ for some elements $t\in T$ and $r_1\in R$. By the assumption, $r_2t=s_1+a$ for some elements $r_2\in R$, $s_1\in S$ and $a\in \ga := \ass_R(T)$. Since $\ass_R(S) = \ga$, $s_2a=0$ for some $s_2\in S$, and so $s_2r_2t= s_2s_1\in S$. Now, $s_2s_1\cdot r= s_2r_2t\cdot r = s_2r_2\cdot tr = s_2r_2\cdot r_1s = s_2r_2r_1\cdot s$. It suffices to take $s'= s_2s_1\in S$ and $r'= s_2r_2r_1\in R$.

(ii) $S\in \Den_l(R, \ga )$: Suppose that $rs=0$ for some elements $ r\in R$ and $s\in S$. Then $r\in \ga =\ass_R(S)$ since $s\in T$. Now, by (i), $S\in \Den_l(R, \ga )$.

(iii) $S^{-1}R\simeq T^{-1}R$ (by Lemma \ref{a7Mar15}.(3)).  $\Box $

Let $\D_1$ be the subalgebra of the Weyl algebra $A_1$ generated by the elements $H$ and $\der$. The algebra $\D_1$ is isomorphic to the skew Laurent polynomial ring  $K[H][\der , \s^{-1}]$ where $\s (H)=H-1$.  Let $A_1^0:=A_1\backslash \{ 0\}$, $A_{1,\der}^0:=A_{1,\der}\backslash \{ 0\}$ and $\D_1^0:=\D_1\backslash \{ 0\}$.

\begin{lemma}\label{a19Mar15}
For each element $a\in \mI_1\backslash F$, there is a natural number $i$ such that $\der^ia\in \D_1^0$.
\end{lemma}

{\it Proof}. By (\ref{acan}), $\der^ia \in \D_1^0+F$ for some $i$. Since $F=\cup_{j\geq 1} \ker (\der^i_{\mI_1}\cdot )$ (where $\der_{\mI_1}^i\cdot : \mI_1\ra \mI_1$, $u\mapsto \der^iu$) and $\der^j\D_1^0\subseteq \D_1^0$, we can enlarge the natural number $i$ such that $\der^ia\in \D_1^0$.
 $\Box $

The next proposition is the key step in finding the rings $Q_{l,cl}(\mI_1)$ and $ Q_{r,cl}(\mI_1)$.

\begin{proposition}\label{1d16Mar15}
Let $\pCC := \pCC_{\mI_1}$, $\pi : \mI_1 \ra \mI_1 / F\simeq A_{1, \der}$, $r\mapsto r+F$, and   $S= \pi^{-1}(A^0_{1, \der})=\mI_1\backslash F$. Then
\begin{enumerate}
\item $S\in \Den_l(\mI_1, F)$ and $S^{-1} \mI_1 \simeq Q(A_1)$.
\item $\D_1^0\in \Den_l(\mI_1, F)$ and ${\D_1^0}^{-1} \mI_1 \simeq S^{-1} \mI_1$.
\item $\pD_1^0:= \pCC_{\mI_1}\cap \D_1^0\in \Den_l(\mI_1, F)$ and ${\pD_1^0}^{-1} \mI_n \simeq{\D_1^0}^{-1} \mI_1$.
\end{enumerate}
Therefore, ${\pD_1^0}^{-1} \mI_1 \simeq{\D_1^0}^{-1} \mI_1\simeq S^{-1} \mI_1 \simeq Q(A_1)$.
\end{proposition}

{\it Proof}. 1. Since  $S_\der:=\{ \der^i\, | \, i\in \N\}\in \Den_l(\mI_1, F)$ and $S_\der \subseteq S$, we have the inclusion $F=\ass (S_\der)\subseteq \ass (S)$. In fact, $F= \ass (S)$ since the algebra $\mI_1/F$ is a domain.  Then $S\in \Den_l(\mI_1, F)$, since $\pi (S) = A_{1,\der}^0\in \Den_l(A_{1,\der}, 0)$, $S_\der \subseteq S$ and $S_\der\in \Den_l(\mI_1, F)$. Now, it is obvious that $S^{-1}\mI_1\simeq \pi (S)^{-1} (\mI_1/F)\simeq (A_{1,\der}^0)^{-1}A_{1,\der}\simeq Q(A_1)$.

2. The inclusion $S_\der\subseteq \D_1^0$ implies that $F=\ass (S_\der ) \subseteq \ass (\D_1^0)$. The factor algebra $\mI_1/F$ is a domain and $\pi|_{\D_1^0}: \D_1^0\ra \D_1^0$ is a bijection, hence $\ass(\D_1^0)\subseteq F$, and so $\ass (\D_1^0)= F$. By Lemma \ref{a19Mar15}, the multiplicative set $\D_1^0$ is dense in $S=\mI_1\backslash F$. By Lemma \ref{c16Mar15}, $\D_1^0\in \Den_l(\mI_1, F)$ and $(\D_1^0)^{-1}\mI_1\simeq S^{-1} \mI_1$.

3.  The inclusion $S_\der \subseteq \pD_1^0$ implies that $F= \ass (S_\der ) \subseteq \ass (\pD_1^0)$. The algebra $\mI_1/ F$ is a domain and $\pi|_{\pD_1^0}: \pD_1^0\ra \pD_1^0$ is a bijection, hence $\ass(\pD_1^0)\subseteq F$, and so $\ass (\pD_1^0)= F$. The right $A_1$-module
 $P_1'$ is a simple one. By \cite{MR-1973-Weyl1}, for every nonzero element $a$ of $A_1$, $\ker (\cdot a_{P_1'})$ is a finite dimensional vector space. In particular, this is the case for all elements $a\in \D_1^0$ (since  $\D_1^0\subseteq A_1^0$). Since $\cap_{i\geq 1} K[\der ]\der^i=0$,   for each  element $a\in \D_1^0$,  we have $\ker(\cdot a_{P_1'})\cap \im(\dot \der^i_{P_1'})=0$ for some $i=i(a)\geq 1$, i.e.  $\ker (\cdot (\der^i a)_{P_1'})=0$. Therefore,  $\der^i a\in \pCC_{\mI_1}\cap \D_1^0= \pD_1^0$ (since $F_{K[\der]}\simeq K[\der ]^{(\N )}$). So, $\pD_1^0$ is dense in $\D_1^0$ and $\ass (\pD_1^0) = \ass (\D_1^0)=F$. By Lemma \ref{c16Mar15}, $\pD_1^0 \in \Den_l(\mI_1, F)$ and $(\pD_1^0)^{-1}\mI_1\simeq (\pD_1^0)^{-1}\mI_1$.
 $\Box $

\begin{corollary}\label{b19Mar15}
For each element $a\in \mI_1\backslash F$ there is a natural number $i$ such that $\der^ia\in \pD_1^0$.
\end{corollary}

{\it Proof}. This was proven in the proof of Theorem \ref{1d16Mar15}.
$\Box $

\begin{theorem}\label{19Mar15}
\begin{enumerate}
\item  $\pQ_{l,cl}(\mI_1) \simeq Q(A_1)$ is a division ring and $\ass (\pCC_{\mI_1}) = F$.
\item $Q_{r,cl}'(\mI_1) \simeq Q(A_1)$ is a division ring and $\ass_r (\pCC_{\mI_1}) = F$.
\end{enumerate}
\end{theorem}

{\it Proof}. 1. By Proposition \ref{1d16Mar15}.(3), $\pD_1^0\in \Den_l(\mI_1, F)$. Since $ \pD_1^0\subseteq \pCC = \pCC_{\mI_1}$ and the algebra $\mI_1/F$ is a domain, we must have $\ass (\pCC ) = F=\ass(\pD_1^0)$. The ideal $F$ is a prime ideal since $\mI_1/F$ is a domain. By Corollary \ref{b19Mar15}, the set $\pbCC :=\pi (\pCC )$ is dense in $\CC_{ \mI_1/F} = A_{1, \der}^0$. By Theorem \ref{28Feb15}.(3), $\pQ_{l, cl} (\mI_1) \simeq
Q_{l, cl}(A_{1, \der})\simeq Q(A_1)$.

2. Applying the involution $*$ to statement 1 and using the fact that $A_1^*\simeq A_1$ we obtain statement 2: $Q_{r,cl}'(\mI_1) = \pQ_{l,cl}(\mI_1)^* \simeq Q(A_1)^*\simeq Q(A_1^*)\simeq Q(A_1)$.
$\Box $

Let $\nabla_1 :=\D_1^*=K[H][\int , \s ]$ where $\s (H)=H-1$ and $\nabla_1^0:= \nabla_1\backslash \{ 0\}$. Applying the involution $*$ of the algebra $\mI_1$ to Proposition \ref{1d16Mar15} and using the facts that $\CC_{\mI_1}'=(\pCC_{\mI_1})^*$, $F^* = F$, $S^* = S$ and $A_1^*\simeq A_1$ (since the Weyl algebra $A_1$ is isomorphic to its dual; in particular $Q(A_1^*)\simeq Q(A_1)$) we obtain the next proposition.

\begin{proposition}\label{a21Mar15}
Let   $S= \mI_1\backslash F$. Then
\begin{enumerate}
\item $S\in \Den_r(\mI_1, F)$ and $ \mI_1S^{-1} \simeq Q(A_1)$.
\item $\nabla_1^0\in \Den_r(\mI_1, F)$ and $ \mI_1{\nabla_1^0}^{-1} \simeq  \mI_1S^{-1}$.
\item $\nabla_1':= \CC_{\mI_1}'\cap \nabla_1^0\in \Den_r(\mI_1, F)$ and $ \mI_n {\nabla_1'}^{-1} \simeq \mI_1{\nabla_1^0}^{-1}$.
\end{enumerate}
Therefore, $ \mI_n {\nabla_1'}^{-1} \simeq \mI_1{\nabla_1^0}^{-1}\simeq  \mI_1S^{-1}\simeq Q(A_1)$.
\end{proposition}

{\bf Descriptions of the sets $\pCC_{\mI_1}$ and $\CC_{\mI_1}'$}.  We are going to give explicit descriptions of the sets  $\pCC_{\mI_1}$ and $\CC_{\mI_1}'$  (Theorem \ref{24Mar15}). They have a sophisticated structure. Let $\G := \{ a=a_0+\sum_{i\geq 1} \int^ia_i+f\, | \, a_0\neq 0, \; {\rm all}\; a_i\in K[H], f\in F\}$. In the proof of Theorem \ref{24Mar15}.(1), it is shown that $\pCC_{\mI_1} \subseteq \cup_{i\geq 1} \der^i\G$ and, for each element $a\in \G$, $\der^ia\in \pCC_{\mI_1}$ for some $i=i(a)\geq 0$. Then map
\begin{equation}\label{GNd}
d: \G \ra \N , \;\; a\mapsto d(a) :=\min \{ i\in \N \, | \, \der^ia\in \pCC_{\mI_1}\}
\end{equation}
is called the {\em left regularity degree function} and the natural number $d(a)$ is called the {\em left regularity degree} of $a$. For each element $a\in \G$, $d(a)$ can be found in finitely many steps, see the proof of Theorem \ref{24Mar15}.(1) where the explicit expression  (\ref{daca})  is given for $d(a)$.
\begin{theorem}\label{24Mar15}

\begin{enumerate}
\item $\pCC_{\mI_1} = \{ \der^{d(a)+i}a\, | \, i\in \N , a\in \G\}$.
\item $\CC_{\mI_1}'= \pCC_{\mI_1}^*$.
\end{enumerate}
\end{theorem}

Before giving a proof of Theorem \ref{24Mar15}, we introduce some definitions.  For each element $a= a_0+\sum_{i\geq 1} \int^ia_i+f\in \G$, the elements $l(a) :=a_0$ and $a_F:= f$  are called the {\em leading term} and the $F$-{\em term} of $a$, respectively. The {\em size} $s(f)$ of the element $f$ is equal to $-1$ if $f=0$, and to $\min \{ m\in \N \, | \, f\in \oplus_{i,j=0}^mKe_{ij}\}$. Then $s(a) := s(a_F)$ is called the {\em size} of $a$. Let $\CR (a_0):= \{ i\in \N \, | \, a_0(i+1)=0\}$, the set of roots of the polynomial $a_0(H+1)$ that are natural numbers. Let $r(a)$ be the maximal element in the set $R(a) := \{ i\in \CR (a_0)\, | \, i>s(a)\}$. If $R(a)=\emptyset$ then $r(a) := \emptyset$.

For each element $a\in \mI_1$, let $\CK_a:= \ker (\cdot a_{P_1'})$.

{\bf Proof of Theorem \ref{24Mar15}}.  Let $\pCC :=\pCC_{\mI_1}$, $P'=P_1'$ and $P'_{\leq i} := \oplus_{j=0}^iK\der^i$ for all $j\geq 0$. Similarly, the vector space $P'_{<i}$ is defined.

1. (i) $\pCC \subseteq \cup_{i\geq 0} \der^i \G$: Let $a\in \pCC$. The element $a$ is a unique sum (\ref{acan}). It suffices to show that there  is $i\leq 0$ such that $a_i\neq 0$ (since then $a= \der^j\g$ where $j=\min \{ i\leq 0\, | \, a_i\neq 0\}$ and $\g \in \G$; this follows from the equalities $\der^k\int^k=1$ and $\der^kF = F$ for all $k\geq 1$). Suppose that $a_i=0$ for all $i\leq 0$, i.e. $a= \sum_{i\geq 1} \int^i a_i+f$ where $f= a_F$. Fix a natural number $n$ such that $n>s(a)$. Then $\cdot a: P_{\leq n} \ra P_{\leq n-1}$, $p\mapsto ap$, and so $\CK_a\neq 0$, a contradiction.

Till the end of the proof let $a=a_0+\sum_{i\geq 1} \int^i a_i+f\in \G$. For all $i\geq s(a)$, $aP_{\leq i}' \subseteq P_{\leq i}'$ and $\der^i\cdot a\equiv \der^i a_0(i+1) \mod P_{\leq i-1}'$. Now, the statement (ii) is obvious.

(ii) $\pG := \G \cap \pCC = \{ a\in \G \, | \, r(a) =\emptyset$ and $ \ker(\cdot a_{P_{\leq s(a)}'})=0\}$. For each element $\pG$, we set $d(a) :=0$. Clearly, $\der^i a\subseteq \pCC$ for all $a\in \pG$ and $i\geq 0$.

Till the end of the proof we assume that $a\in \G\backslash \pG$. By the statement (ii), there are two cases:

(a) $r(a) \neq \emptyset $ (i.e. $R(a) \neq \emptyset$), and

(b) $r(a) = \emptyset$ (i.e. $R(a) =\emptyset$) and $\ker (\cdot a_{P_{\leq s(a)}}')\neq 0$.

In the case (a), $\CK_a \subseteq P'_{\leq r(a)}$, and so $\der^is\in \pCC$ for all $i>r(a)$. In the case (b), $\CK_a \subseteq P'_{\leq s(a)}$, and so $\der^is\in \pCC$ for all $i>s(a)$. This proves that the function $d$ is well-defined (see (\ref{GNd})) and that $\pCC_{\mI_1} = \{ \der^{d(a)+i}a\, | \, i\in \N , a\in \G\}$, in view of  the statements (i) and (ii).  Then
\begin{equation}\label{daca}
d(a) = \begin{cases}
\min \{ i\, | \, 0\leq i\leq r(a)+1, \; \CK_a\cap(\bigoplus_{j=i}^{r(a)+1}K\der^j)=0\}& \text{in the case (a)},\\
\min \{ i\, | \, 0\leq i\leq s(a)+1, \; \CK_a\cap(\bigoplus_{j=i}^{s(a)+1}K\der^j)=0\}& \text{in the case (b)}.
\end{cases}
\end{equation}
2. Statement 2 is obvious. $\Box$





\small{

Department of Pure Mathematics

University of Sheffield

Hicks Building

Sheffield S3 7RH

UK

email: v.bavula@sheffield.ac.uk}

\end{document}